%
\magnification\magstep1

\def\forces{\parallel\!\!\! -}


\def\hexnumber#1{\ifcase#1 0\or1\or2\or3\or4\or5\or6\or7\or8\or9\or
	A\or B\or C\or D\or E\or F\fi }

\font\teneuf=eufm10
\font\seveneuf=eufm7
\font\fiveeuf=eufm5
\newfam\euffam
\textfont\euffam=\teneuf
\scriptfont\euffam=\seveneuf
\scriptscriptfont\euffam=\fiveeuf
\def\frak{\fam\euffam \teneuf}

\font\tenmsx=msam10
\font\sevenmsx=msam7
\font\fivemsx=msam5
\font\tenmsy=msbm10
\font\sevenmsy=msbm7
\font\fivemsy=msbm5
\newfam\msxfam
\newfam\msyfam
\textfont\msxfam=\tenmsx  \scriptfont\msxfam=\sevenmsx
  \scriptscriptfont\msxfam=\fivemsx
\textfont\msyfam=\tenmsy  \scriptfont\msyfam=\sevenmsy
  \scriptscriptfont\msyfam=\fivemsy
\edef\msx{\hexnumber\msxfam}

\mathchardef\upharpoonright="0\msx16
\let\restriction=\upharpoonright

\def\restrict{{\restriction}}

\def\qed{{\vcenter{\hrule height.4pt \hbox{\vrule width.4pt height5pt
 \kern5pt \vrule width.4pt} \hrule height.4pt}}}
\def\notin{{\in}\kern-5.5pt / \kern1pt}
\def\ok{\vbox{\hrule height 8pt width 8pt depth -7.4pt
    \hbox{\vrule width 0.6pt height 7.4pt \kern 7.4pt \vrule width 0.6pt height 7.4pt}
    \hrule height 0.6pt width 8pt}}
\def\nt{{\leq}\kern-1.5pt \vrule height 6.5pt width.8pt depth-0.5pt \kern 1pt}
\def\sd{{\times}\kern-2pt \vrule height 5pt width.6pt depth0pt \kern1pt}
\def\zp#1{{\hochss Y}\kern-3pt$_{#1}$\kern-1pt}

\def\extend { \hat {\; \; } }
\def\sm{{\smallskip}}

\def\no{{\noindent}}
\def\la{{\langle}}

\def\sub{\subseteq}
\def\alm{\sub ^*}

\def \o {\omega }
\def \fun {{^\omega \omega }}

\font\capit=cmcsc10 scaled\magstep0

\overfullrule=0pt
\openup1.5\jot

\def \inf {[\o ]^\o }
\def \fin {{^{<\omega }\omega }}
\def \o {\omega } \def \sub {\subseteq }
\def \fun {{^{\omega }\omega }}
\def \k {\kappa } 
\def \l {\langle } \def \r {\rangle }
\def \g {\gamma } \def \s {\sigma } \def\t {\tau } 
\font\gross=cmbx10 scaled \magstep1
\font\sgross=cmbx10 scaled \magstep2
  
\def \Y {{\cal Y}} \def \P {{\cal P}} \def \A {{\cal A}} \def \B
{{\cal B}}  \def\D{{\cal D}} 
\def\square{\qed }

\def\Vdash{\forces }

\def \fin {{^{<\omega }\omega }}
\def \o {\omega } \def \sub {\subseteq }
\def \alm {\subseteq ^*}
\def \fun {{^{\omega }\omega }}
\def \inf {[\o ]^\o }
\def \l {\langle } \def \r {\rangle }
\def \s {\sigma }
\def \k {\kappa}
\def \a {\alpha }
\def \b {\beta }
\def \d {\delta }
\def \t {\tau }
\font\hgross=cmbx10 scaled \magstephalf
\font\gross=cmbx10 scaled \magstep1
\font\sgross=cmbx10 scaled \magstep2
\def \iter {( P_\a , \dot Q_\b : \a \leq \o _2 ,\b < \o _2 )}  
\def \P {{\cal P}} \def \Y {{\cal Y}}
\def \dom {\hbox{dom}}
\def \h {{\frak h}}

\newcount\skewfactor
\def\mathunderaccent#1#2 {\let\theaccent#1\skewfactor#2
\mathpalette\putaccentunder}
\def\putaccentunder#1#2{\oalign{$#1#2$\crcr\hidewidth
\vbox to.2ex{\hbox{$#1\skew\skewfactor\theaccent{}$}\vss}\hidewidth}}
\def\name{\mathunderaccent\tilde-3 }
\def\Name{\mathunderaccent\widetilde-3 }



\def\hexnumber#1{\ifcase#1 0\or1\or2\or3\or4\or5\or6\or7\or8\or9\or
	A\or B\or C\or D\or E\or F\fi }

\font\teneuf=eufm10
\font\seveneuf=eufm7
\font\fiveeuf=eufm5
\newfam\euffam
\textfont\euffam=\teneuf
\scriptfont\euffam=\seveneuf
\scriptscriptfont\euffam=\fiveeuf
\def\frak{\fam\euffam \teneuf}

\font\tenmsx=msam10
\font\sevenmsx=msam7
\font\fivemsx=msam5
\font\tenmsy=msbm10
\font\sevenmsy=msbm7
\font\fivemsy=msbm5
\newfam\msxfam
\newfam\msyfam
\textfont\msxfam=\tenmsx  \scriptfont\msxfam=\sevenmsx
  \scriptscriptfont\msxfam=\fivemsx
\textfont\msyfam=\tenmsy  \scriptfont\msyfam=\sevenmsy
  \scriptscriptfont\msyfam=\fivemsy
\edef\msx{\hexnumber\msxfam}

\mathchardef\upharpoonright="0\msx16
\let\restriction=\upharpoonright

\def\extend { \hat {\; \; } }
\def\restrict{{\restriction}}

\def\qed{{\vcenter{\hrule height.4pt \hbox{\vrule width.4pt height5pt
 \kern5pt \vrule width.4pt} \hrule height.4pt}}}
\def\notin{{\in}\kern-5.5pt / \kern1pt}
\def\ok{\vbox{\hrule height 8pt width 8pt depth -7.4pt}}
    
\def\square{\qed }

\def\Vdash{\forces }

\def \fin {{^{<\omega }\omega }}
\def \o {\omega } \def \sub {\subseteq }
\def \alm {\subseteq ^*}
\def \fun {{^{\omega }\omega }}
\def \inf {[\o ]^\o }
\def \l {\langle } 
\def \r {\rangle }
\def \s {\sigma }
\def \a {\alpha }
\def \b {\beta }
\def \d {\delta }
\font\hgross=cmbx10 scaled \magstephalf
\font\gross=cmbx10 scaled \magstep1
\font\sgross=cmbx10 scaled \magstep2
\def \iter {\l P_\a , \name{ Q_\b} : \a \leq \o _2 ,\b < \o _2 \r}  
\def \P {{\cal P}} \def \Y {{\cal Y}}
\def \dom {\hbox{dom}}
\def \h {{\frak h}}
\def \sm {\smallskip}
\def\no{\noindent}
\def\B{{\cal B}}
\def\E{{\name{ E}}}
\def\g{{\gamma}}
\def\D{{\name{D}}}

\noindent {\sgross The distributivity numbers of of $\P (\o )$/fin and
its square}

\bigskip \bigskip \bigskip

\noindent Saharon Shelah\footnote{$^1$}{The author is supported by the
Basic Research Foundation of the Israel Academy of Sciences;
publication 494.}

\smallskip

\item{}{Institute of Mathematics, Hebrew University, Givat Ram, 91904
Jerusalem, ISRAEL}

\smallskip

\noindent Otmar Spinas\footnote{$^2$}{The author is supported by the
Swiss National Funds.}

\smallskip

Mathematik, ETH-Zentrum, 8092 Z\"urich, SWITZERLAND

\smallskip

\item{}{Institut f\"ur Mathematik, Humboldt-Universit\"at, Unter den
Linden 6, 10099 Berlin, GERMANY ({\it current address})}

\bigskip \bigskip

{\narrower ABSTRACT: We show that in a model obtained by forcing with
a countable support iteration of Mathias forcing of length $\o _2$,
the distributivity number of ${\cal P}(\o )$/fin is $\o _2$, whereas
the distributivity number of r.o.$({\cal P}(\o )$/fin)$^2$ is $\omega
_1$.  This answers an old problem of Balcar, Pelant and Simon, and
others.

}

\bigskip \bigskip

\def \sub {\subseteq}
\def \o {\omega }

\def \A {{\cal A}}
\def \d {\delta } \def \t {\tau }
\def \tree {{^{<\o }2}}
\def \k {\kappa }
\def \b {\beta }
\def \a {\alpha }
\def \P {{\cal P}(\omega )}
\def \h {{\frak h}}
\def \ha {{\frak h}(2)}
\def \la {\lambda }

{\gross Introduction}

\bigskip 

A complete Boolean algebra $( B,\leq ) $ is called $\k
$-distributive, where $\k $ is a cardinal, if and only if for every
family $\l u_{\a i}: i\in I_\a , \a <\k \r $ of members of $B$ 

$$\prod _{\a < \k }\sum _{i\in I_\a } u_{\a i} = \sum _{f\in \prod _{\a
<\k } I_\a } \prod _{\a <\k } u_{\a f(\a )} $$

\noindent holds. It is well-known (see [J, p.152]) that every partially
ordered set $( P,\leq ) $ which is separative can be
densely embedded in a unique complete Boolean algebra, which is
usually denoted with r.o.$(P)$. The distributivity number of $(
P,\leq) $ is the defined as the least $\kappa $ such that r.o.$(P)$
is not $\kappa $-distributive. It is well-known (see [J, p.158])
that the following four statements are equivalent:

\smallskip

\item{(1)}{r.o.$(P)$ is $\k $-distributive.}

\item{(2)}{The intersection of $\k $ open dense sets in $P$ is dense.}

\item{(3)}{Every family of $\k $ maximal antichains of $P$ has a
refinement.} 

\item{(4)}{Forcing with $P$ does not add a new subset of $\k $.}

\smallskip

The distributivity number of the Boolean algebra $\P $/fin is denoted with
${\h }$. This
cardinal was introduced in [BPS], where it has been shown that $\o
_1 \leq \h \leq 2^\o $  and the axioms of ZFC do not decide
where exactly $\h $ sits in this interval. 

For $\la $ a cardinal let $\h (\la )$ be the distributivity number of
$(\P $/fin$)^\la $, where by 
$(\P $/fin$)^\la $ we mean the full $\la $-product of $\P $/fin in the
forcing sense. That is, $p\in (\P $/fin$)^\la $ if and only if $p: \la
\rightarrow \P $/fin $\setminus \{ 0\} $. The ordering is
coordinatewise. 

Trivially, $\h (\la )\geq \h (\gamma )$ holds whenever $\la <
\gamma $. In fact, if $\l D_\a :\a<\h (\lambda )\r$ is a family of dense
open subsets of $(\P $/fin$)^\lambda $ whose intersection is not dense,
then, letting $D_\a'=\{ p\in (\P \hbox{/fin})^\gamma :p\restrict \la
\in D_\a\}$, clearly the $D_\a'$ are dense open in $(\P $/fin$)^\gamma
$ and their intersection is not dense.

Since $\h\leq 2^\o$, this implies that under CH the sequence $\l \h
(\la):\la \in \hbox{{\bf Card}}\r$ is constant with value $\aleph _1$.
In [BPS, 4.14(2)] we read: ``We do not know of any further properties
of this sequence.'' The most elementary question which arises, and
which was explicitly asked by several people, is whether consistently
this sequence is not constant. In this paper we give a positive answer
by proving the consistency of $\h (2)<\h$ with ZFC. In a sequel paper,
for every $n<\o$ we will construct a model for $\h(n+1)<\h(n)$. In all
these models the continuum will be $\aleph _2$, and hence the above
sequence will be two-valued. The question of whether more values are
possible is tied up with the well-known problem of how to make the
continuum bigger than $\aleph _2$, not using finite-support forcing
iterations.

The natural forcing to increase $\h $ is Mathias forcing. We will show
that in a model obtained by forcing with a countable support iteration
of length $\o _2 $ of Mathias forcing over a model for CH, $\ha $
remains $\o _1 $. 

There exists an equivalent game-theoretic definition of $\h (\la )$,
which we will use in the sequel. For any ordinal $\a $ and any partial
ordering $P$ let us consider the following game $G(P,\a )$ of length
$\a $: Player I and II alternately choose elements $p^I_\b ,
p^{II}_\b \in P$, $\b < \a $, such that for $\b <\b ' < \a $: $p^I_\b
\geq p^{II}_\b
\geq p^I_{\b '} \geq p^{II}_{\b '}$. In the end, player II wins if and
only if the sequence of moves has no lower bound (this might happen if
at some step $\b<\a$, player I does not have a legal move).

We claim that $\h (\la ) $ is the minimal cardinal $\k $ such that in
the game $G((\P $/fin$)^\la , \k )$, player II has a winning strategy.
For one direction, suppose we are given dense open sets $\l D_\a :\a
<\k \r $ in $(\P $/fin$)^\la $ such that $D= \bigcap \{ D_\a :\a <\k
\} $ is not dense. By the homogeneity of $(\P $/fin$)^\la $ we may
assume that $D$ is empty. In fact, if $D$ contains no extension of
$p$, choose $\l f_\a:\a<\la\r$ such that $f_\a:p(\a )\rightarrow \o$
is one-one and onto. Replace $D_\a$ by $D_\a'=\{ \l f_\a[q(\a)]:\a<\la\r:
q\in D_\a\hbox{ and }q\leq p\}$. Then the $D_\a'$ are open dense and
their intersection is empty. Now define a strategy for II in $G((\P
$/fin$)^\la , \k )$ as follows: In his $\a $th move let II play
$p^{II}_\a
\in D_\a $ such that $p^{II}_\a \leq p^{I}_\a $. This is clearly a
winning strategy.

Conversely, let $\s $ be a winning strategy for II in $G((\P
$/fin$)^\la , \k )$. We will make use of (3) above. We define maximal
antichains $\l \A _\a : \a <
\gamma \leq \k \r $ in $(\P $/fin$)^\la $ such that if
$\a <\beta <\gamma $, then $\A _\beta $ refines $\A _\a $, and for every
$p_\beta \in \A _\beta $, if $p_\a \in \A _\a $ is the unique member
with $p_\a \geq p_\beta $, then $\l p_\a :\a \leq \beta \r $ are
responses by $\s $ in an initial segment of a play consistent with $\s
$. Suppose $\l \A _\a
: \a <\delta \r $ has been constructed and $\delta <\k $ is a limit.
If this sequence has no refinement we are done, otherwise let $\B $ be
one. Now it is easy to construct $\A _\delta $ as desired, namely
consisting of responses by $\s $ to plays of length $\delta +1$ with
last coordinate an extension of a member of $\B $. If $\delta $ is a
successor, construct $\A _\delta $ similarly, where now $\B =\A
_{\delta -1}$. It is clear that
this construction stops at some $\gamma \leq \k $, as otherwise we
could find a play consistent with $\s $ in which II loses.




\bigskip

{\gross 1. Mathias forcing and Ramsey ultrafilters}

\bigskip

Conditions of Mathias forcing are pairs $( u,a) \in [\o ]^{<\o }
\times [\o ]^\o $ such that max$(u) < $ min$(a)$. The ordering is
defined as follows: $( u,a) \leq ( v,b) $ if and only if $v\sub u \sub
v\cup b$ and $a\sub b$. Mathias forcing will be denoted by $Q$ in this
paper. Given $p\in Q$ we will write $p=(u^p,a^p)$. 

If ${ D}$ is a filter on $\o $ containing no finite sets, then
$Q({ D})$ denotes Mathias forcing relativized to ${ D}$, that
is, $(u,a)\in Q({ D})$ iff $(u,a)\in Q$ and $a\in { D}$, and
the order is as for $Q$. Note that any two conditions in $Q({ D})$
with the same first coordinate are compatible. Therefore, $Q({
D})$ is $\sigma $-centered, that is, a countable union of centered
subsets. It is well-known that Mathias forcing can be decomposed as
$Q=Q'\ast \name{ Q''}$, such that $Q'$ is $\P $/fin
and $\name{ Q''}=Q(\name{ G'})$, where $\name{ G'}$ is a name for the generic
filter added by $\P$/fin. In fact, since $Q'$ is $\sigma$-closed and
hence does not add reals, the map sending $(u,a)$ to $(a,(u,a))$ is a
dense embedding of $Q$ in $Q'\ast \name{ Q''}$.
The generic filter for $\name{ Q''}$, which
determines the Mathias real, will be denoted $\name{ G''}$. 
Here and in the sequel we do not distinguish between a member of $\P
$/fin and its representatives in $\P$.  {\bf The above notation will
be used throughout the paper.}

\sm 

The Rudin-Keisler order $\leq _{RK}$ for ultrafilters on $\o $
is defined by: $D\leq _{RK} U $ iff there exists a function
$f:\o \rightarrow \o $ such that $D=\{ X\sub \o :f^{-1}[X] \in U\} $.
In this case $D$ is called a projection of $U$ and it is denoted by
$f_*(U)$. If $D\leq _{RK}U$ and $U\leq _{RK}D$, we call $U$ and $D$
RK-equivalent. By a result of M.E. Rudin (see [R] or [J, 38.2.,
p.480]), in this case there
exists a bijection $f:\o \rightarrow \o $ such that $D=f_*(U)$. Then
we say that $D$ and $U$ are RK-equivalent by $f$. 

A nonprincipal ultrafilter $D$ on $\o $ is called a Ramsey ultrafilter iff for
every $n,k<\o$ and every partition $F:[\o]^n\rightarrow k$ there
exists $H\in D$ homogeneous for $F$, that is, $F\restrict [H]^n$ is
constant. An equivalent definition is as follows (see [J, p.478]): 
$D$ as above is Ramsey
iff for every partition of $\o$ into pieces not in the filter there
exists a filter set which meets each piece at most once. Clearly such
a filter is a $p$-point, that is, for every countable subset of the
filter there exists a filter set which is almost contained in every
member of it. 

We will use yet another equivalent definition of Ramsey ultrafilter.
Let $D$ be a nonpricipal ultrafilter.
A function $f\in \fun $ is called
unbounded modulo $D$ if $\{ n:f(n)>k\} \in D$ for every $k<\o $;
moreover $f$ is called one-to-one modulo $D$ if its restriction to some
member of $D$ is one-to-one. Then $D$ is a Ramsey ultrafilter
iff every function unbounded modulo $D$ is one-to-one modulo $D$ (see
[J, 38.1.,p.479]).

In the following lemma, a forcing $P$ is called $\fun$-bounding iff every
function in $\fun$ in the extension $V^P$ is bounded by some function
in $V$. Moreover, an ultrafilter $D$ in $V$ is said
to generate an ultrafilter in $V^P$ iff the collection of subsets of
$\o$ which belong to $V^P$ and contain an element of $D$ is an
ultrafilter in $V^P$.

\sm

{\bf Lemma 1.1.}  {\it Suppose $D_1 ,D_2$ are Ramsey ultrafilters
which are not RK-equivalent. Let $P$ be a proper, $\fun $-bounding forcing
such that for every filter $G\sub P$ which is $P$-generic over $V$,
$D_1$ and $D_2$ generate ultrafilters in $V[G]$. Then in $V[G]$, $D_1$
and $D_2$ generate Ramsey ultrafilters which are not RK-equivalent.}

\sm

{\it Proof:}  Firstly, we show that $D_1, D_2$ are Ramsey ultrafilters  
in $V[G]$. Here and in the sequel, we denote the ultrafilters generated
by $D_1,D_2$ in $V[G]$ by $D_1,D_2$ as well. By properness, every
$X\in [V]^\o \cap V[G]$ is covered by a countable set in $V$. Hence
$D_1,D_2$ generate $p$-points in $V[G]$. In $V[G]$, let $\l a_n:n<\o
\r $ be a partition of $\o $ such that $a_n \not\in D_1 $, for all
$n<\o $. As $D_1$ is a $p$-point, there exists $X\in D_1\cap V$ such
that $|X\cap a_n|<\o $, for all $n<\o $. Let $f\in \fun $ be defined
by: $f(n+1)>f(n)$ is minimal such that every $a_k$ with $a_k\cap
f(n)\ne \emptyset $ satisfies $a_k \cap (X \setminus f(n+1))=\emptyset
$. As $P$ is $\fun $-bounding, we may find a strictly increasing $g\in
\fun \cap V$ such that for every $n<\o $, $[g(n),g(n+1))\cap $
range$(f)$ has at least one element. $D_1$ contains exactly one of the
three sets $\bigcup \{ [g(3n+i),g(3n+i+1)):n<\o \} $, where $i\in \{ 0,1,2\}$.
We denote this set by $Y$. 
As $D_1$ is Ramsey in $V$, there exists $Z\in
D_1\cap V$ such that $Z\sub X\cap Y$ and
$|[g(n),g(n+1)) \cap Z|\leq 1$, for all $n<\o $.  We have to verify
that $|Z\cap a_n|\leq 1$, for every $n$. Let $k,l\in Z\cap a_n$. Then
$k,l\in X\cap a_n$. By construction of $f$, there is $n_1$ such that
$X\cap a_n\sub [f(n_1),f(n_1 +2))$. By construction of $g$ and since
$f$ is increasing, there is $n_2$ such that $f(n_1),f(n_1
+1),f(n_1+2)\in [g(n_2),g(n_2+3))$. By construction of $Z$, there is
$n_3\in\{ n_2,n_2+1,n_2+2\}$ such that $k,l\in [g(n_3),g(n_3+1))$.
Since $|[g(n_3),g(n_3+1)) \cap Z|\leq 1$, we have that $k=l$.

Secondly, we show that $D_1,D_2$ do not become RK-equivalent in
$V[G]$. Otherwise, in $V[G]$ we had 
a bijection $f:\o\rightarrow \o$ such that $f_*(D_1)=D_2$. Let
$f_1\in \fun $ be defined such that $f_1(n+1)>f_1(n)$ is minimal with

$$f_1(n+1)\geq \max [\{ f(k):k<f_1(n)\} \cup \{ f^{-1}(k):k<f_1(n)\}
]\, .$$

As $P$ is $\fun $-bounding, we may find a strictly increasing $g\in
\fun \cap V$ such that for every $n<\o $, $[g(n), g(n+1))\cap \hbox{
range} (f_1)$ has at least two elements. Each of $D_1$ and $D_2$ contains
one of the three sets 

$$C_i=\bigcup \{ [g(3n+i), g(3n+i+1 )): n<\o \} \, ,$$

\noindent where $i\in \{ 0,1,2\} $. Suppose $C_i\in D_1$ and $C_j \in
D_2 $. By Ramseyness in
$V$, there exist $X\in D_1\cap V$, $Y\in D_2 \cap V$ such that $X\sub
C_i$, $Y\sub C_j$ and $|X\cap [g(3n+i), g(3n+i+1))|\leq 1 $,
$|Y\cap [g(3n+j), g(3n+j+1))|\leq 1$, for all $n<\o $.
Let $x_n$ be the unique element of $X\cap [g(3n+i),g(3n+i+1))$ in the
case that this set is not empty, and let $y_n$ be the unique element of
$Y\cap [g(3n+i-1), g(3n+i+1))$ if this set is not empty. 
Note that by
construction, $f(x_n)\in [g(3n+i-1), g(3n+i+1))$. Hence $\{ x_n:
f(x_n)=y_n\} \in D_1 $, as otherwise $f$ would map a set in $D_1$ to a
set disjoint to a member of $D_2$. Consequently, $\{
y_n:f(x_n)=y_n\}\in D_2$. Choose $X_1\in D_1\cap V$ and $Y_1\in
D_2\cap V$ such that
$X_1\sub \{ x_n:f(x_n)=y_n\} $ and $Y_1\sub\{ y_n:f(x_n)=y_n\} $.
Define

$$f'=\{ (x,y):\exists n(x\in [g(3n+i),g(3n+i+1))\cap X_1\wedge y\in
[g(3n+i-1),g(3n+i+2))\cap Y_1)\}\, $$

\no Then $f'\in V$ and $f'$ is a map with dom$(f')=f^{-1}[f[X_1]\cap
Y_1]\in D_1$ and $f'(x)=f(x)$ for all $x\in$ dom$(f')$. Therefore,
$f'$ witnesses in $V$ that $D_1,D_2$ are RK-equivalent, a
contradiction. 

\hfill $\qed$

\sm

In the sequel we will have the
following situation: Given are two models of ZFC, $V_0\sub V_1$, and
in $V_1$ we have $D$ which is an ultrafilter on $(\inf )^{V_0}$. That is,
$D\sub (\inf )^{V_0}$ is a filter and for every $a\in (\inf )^{V_0}$, either
$a\in D$ or $\o \setminus a\in D$. Then we call $D$ Ramsey if every
function in $V_0$ which is unbounded modulo $D$ is one-to-one modulo
$D$. We will say that some real $r\in (\inf) ^{V_1}$ induces $D$ if
$D=\{ a\in (\inf)^{V_0}:r\sub ^*a\}$.

An easy genericity argument together with the $\sigma$-closedness of
$\P$/fin shows that $\forces _{\P \hbox{/fin}} \,\name{ G'}$ is a Ramsey
ultrafilter.

In [M], Mathias has shown that $r\in \inf$ is Mathias generic over $V$
if and only if $r$ is an almost intersection of a 
${\cal P}(\o)$/fin-generic filter $G'$, that is, $r\sub ^*a$ for all
$a\in G'$. It follows that every infinite subset of a Mathias generic
real is Mathias generic as well. This will be used in the proof of the
following well-known fact.

\sm

{\bf Lemma 1.2.} {\it Let $(N,\in)$ be a countable model of ZF$^-$ (in
particular, $N$ must be able to prove the above mentioned result of
Mathias). If $p\in Q\cap N$ there exists $q\in Q$ such that $q\leq p$,
$u^p=u^q $, and for every $a\in \inf$ with $u^q\sub a\sub u^q\cup
a^q$, $a$ is Mathias generic over $N$. In particular, $q$ is
$(N,Q)$-generic below $p$.}

\sm

{\it Proof:} Since $N$ is countable, in $V$ we may find $b\in \inf$
which is Mathias generic over $N$ and contains $p$ in its induced
generic filter, that is, $u^p\sub b\sub u^p\cup a^p$. Let $q=(u^p,
b\setminus u^p)$. Then every $a$ as in the Lemma is an infinite subset of
$b$, and hence Mathias generic over $N$.\hfill $\qed$

\bigskip

{\gross 2. Outline of the proof}

\bigskip

Let $V$ be a model of CH and let $\l P_\a ,\name{ Q_\b} : \a \leq \o
_2,\b < \o _2\r $ be a countable support iteration of Mathias forcing,
that is $\forall \a < \o _2 $, $\Vdash _{P_\a } $ ``$\name{ Q_\a} $ is
Mathias forcing''. {\bf This notation will be kept throughout the paper.}

The following theorem is folklore. In the proof, a set $C\sub
\o _2$ will be called $\o _1$-club if $C$ is unbounded in $\o _2$
and closed under
increasing sequences of length $\o _1$.

\sm
 
{\bf Theorem 2.1.}  {\it If $G$ is $P_{\o _2}$-generic over $V$, where
$V\models CH$, then $V[G]\models \h =\o _2$.}

\sm

{\it Proof:}  In $V[G]$ let $\l D_\nu :\nu < \o _1 \r $ be a family of
open dense subsets of $\P $/fin $\setminus \{ 0\} $. By a standard
L\"owenheim--Skolem argument,
for every $\a $ belonging to some $\o _1$-club $C\sub \o _2$, for every $\nu
<\o _1$ it is true that $D_\nu \cap V[G_\a ]$ belongs to $V[G_\a ]$
and is open dense in $(\P $/fin)$^{V[G_\a ]} \setminus \{ 0\} $. Now
for given $A\in
(\P $/fin$)^{V[G]}\setminus \{ 0\} $, by properness and genericity
there exists $\a \in C$ such that $A\in G(\a )'$, where $G(\a )$ is
the $\name{ Q_\a} [G_\a ]$-generic filter determined by $G$ and $G(\a )'$
is its first component according to the decomposition of Mathias
forcing defined in $\S 1$. As $\a \in C$, $G(\a
)'$ 
clearly meets every $D_\nu $, $\nu <\o _1$. But now $r_\a $, the $\name{
Q_\a }$-generic real (determined by $G(\a )''$) is below each member of
$G(\a )'$, hence below $A$ and in $\bigcap _{\nu <\o _1}D_\nu $. This
proves that $\bigcap _{\nu <\o _1}D_\nu $ is dense. \hfill $\qed $

\sm

All the rest of this paper is to prove:

\sm

{\bf Theorem 2.2.}  {\it In the notation of Theorem 2.1, $V[G]\models \ha
=\o _1 $.}

\sm

The proof consists of the following two propositions.  By $S^2_1$ we
will denote the ordinals in $\o _2$ of cofinality $\o _1$. We will
tacitly use the well-known results from [B, $\S 5$], where it has been
shown that for $\a <\o _2$ we can define a quotient forcing $P_{\o
_2}/\name{ G_\a} $, also denoted $P_{\a \o _2}$, where $\name{ G_\a} $ is a
$P_\a $-name for the $P_\a$-generic filter.

\sm

{\bf Proposition 2.3.} {\it There exists an $\o _1$-club $C\sub S^2_1$
such that for every $\a \in C$ the following holds: If $\name{ r} $ is a
$P_{\o _2}/\name{ G_\a} $-name such that $\Vdash _{P_{\o _2}/\name{ G_\a
}}\hbox{ ``}
\name{ r} $ induces a Ramsey ultrafilter on $([\o ]^\o )^{V[\name{ G_\a} ]} $'',
then $\Vdash _{P_{\o _2}/\name{ G_\a }}\, \name{ r} \in V[\name{ G_{\a
+1}}]$.}

\sm

{\bf Proposition 2.4.} {\it Suppose that $V\models CH$ and $\name{ r} $ is a
$Q$-name such that $\Vdash _Q$ ``$\name{ r}
$ induces a Ramsey ultrafilter $\D $ on $([\o ]^\o )^{V}$''. 
Then $\Vdash _Q\hbox{
``}  \D $ and $ \name{ G'}$
are RK-equivalent by some function $f \in ({^\o }\o )\cap V\hbox{''.}$}

\sm

It is easy to see that Theorem 2.2 follows from Propositions 2.3 and
2.4: Fix
$C$ as in Proposition 2.3. In
$V[G]$ define a winning strategy for player II in the game
$G((\P$/fin$)^2, \o _1)$ as follows: 

\smallskip

{\narrower{\no Play in such a way that whenever $\l (p^I_\nu ,p_\nu ^{II}) :
\nu <\o _1 \r $ is a play, there exists $\a \in C$ such that $\l
p^{II}_\nu (0) :\nu <\o _1 \r $ and $\l p^{II}_\nu (1):\nu <\o _1 \r $
generate Ramsey ultrafilters on $([\o ]^\o )^{V[G_\a ]}$ which are not
RK-equivalent by any $f \in ({^\o }\o )^{V[G_\a ]} $.

}}

\smallskip

First we show that such a strategy exists in $V[G]$. Then we show that
it is winning. We work in $V[G]$. For $x\in V[G]$, let $o(x)=\min\{
\a<\o_2 :x\in V[G_\a ]\}$. Let $\Gamma : \o _1\rightarrow (\o _1)^2$ be a
bijection such that $\Gamma (\a )=(\b ,\delta )$ implies $\b\leq\a$.
For each $\a <\o _2$, $V[G_\a]\models $ CH. Hence we can choose $g_\a
:\o _1\rightarrow V[G_\a]$ which enumerates all triples $(a,\pi ,f)\in
V[G_\a]$ such that $a\in \inf$, $\pi :[\o]^n\rightarrow k$ for some
$n,k<\o$, and $f\in\fun$. In his $\a$th move, II plays $(p^{II}_\a
(0),p^{II}_\a (1))\leq (p^{I}_\a (0),p^{I}_\a (1))$ such that, {\it
if} $\Gamma (\a)=(\b ,\d)$, $\xi\in C$ is minimal with $\xi \geq
\sup\{ o((p^{I}_\nu
(0),p^{I}_\nu (1))):\nu <\b\} $, and $(a,\pi ,f)=g_\xi (\delta )$,
{\it then} for $i\in \{ 0,1\}$ we have:

\item{(1)} $p^{II}_\a(i)\sub a$ or $p^{II}_\a(i)\cap a=\emptyset $,

\item{(2)} $p^{II}_\a(i)$ is homogeneous for $\pi$,

\item{(3)} $f[p^{II}_\a(0)]\cap p^{II}_\a(1)=\emptyset$.

\no As $C$ is $\o _1 $-club, it is easy to verify that 
this strategy is as desired.

Suppose that $\l p_\nu :\nu <\o_1\r$ are moves of player II which are
consistent with this strategy. Suppose this play is won by I. Hence
there exists $(r_0,r_1)\in (\inf )^2\cap V[G]$ with $(r_0,r_1)\leq
p_\nu$, for all $\nu <\o_1$. So we get $\a \in C$, and Ramsey
ultrafilters $G_i$ on $(\inf )^{V[G_\a]}$, for $i<2$, such that $G_i$
is generated by $\l p_\nu (i) :\nu <\o_1\r$, and $G_0$ is not
RK-equivalent to $G_1$ by any $f\in \fun \cap V[G_\a]$.  Then $G_i$ is
generated by $r_i$. By Proposition 2.3 we obtain that $r_i$ belong to
$V[G_{\a +1}]$, and hence by Proposition 2.4, $G_0$ and $G_1$ are both
RK-equivalent to $G(\a )'$ by some $f\in \fun\cap V[G_\a]$. By
construction this is impossible.  By the game-theoretic
characterization of $\h (2)$ (see Introduction), this implies
$V[G]\models \ha =\o _1$.

\bigskip

{\gross 3. Iteration of Mathias forcing}

\bigskip

Throughout this section $\l P_\a ,\name{ Q_\b} :\a\leq \gamma
,\b<\gamma\r$ denotes a countable support iteration of Mathias forcing
of length $\g$. By [Sh{\it b}, p.96ff.] we may assume that elements of
$P_\g $ are hereditarily countable. {\bf We shall always assume this
in the sequel.} For $p\in P_\g$, the collection of $\beta\in \g$ such
that in the transitive closure of $p$ there exists a $P_\beta$-name
for a condition in $\name{ Q_\b}$, is denoted with cl$(p)$. By our
assumption, cl$(p)$ is a countable subset of $\g$. Note that if $\l
r_\a:\a<\g\r$ is a sequence of $P_\g$-generic Mathias reals, then only
$\l r_\a:\a\in $ cl$(p)\r$ are needed in order to evaluate $p$.  
Letting $a^*=$ cl$(p)$, we can define $P_{a^*}$ as the countable
support iteration of Mathias forcing with domain $a^*$. So $P_{a^*} $
is isomorphic to $P_\delta$, where $\delta =$ o.t.$(a^*)$. The
question arises whether we can view $p$ as a condition in $P_{a^*}$.
It should be noticed that this is not trivially the case.

In this section we prove that $P_\gamma$ has a dense subset $P_\gamma
'$ which can be equipped with an order $\leq '$, such that forcing
with $(P_\gamma ,\leq )$ is equivalent to forcing with $(P_\g ',\leq
')$, and the definition of $(P_\g
',\leq ')$ is absolute for $\Pi ^1_1$-correct models of
ZF$^-$ (up to some trivial restrictions). This will be used in the
following sections to show that potential counterexamples to
Propositions 2.3 and 2.4 must be added by an iteration of countable
length (see Lemma 4.2). In particular, it will be obvious that if $p\in
P_\gamma '$, then $p\in P_{a^*}'$, where $a^*=$ cl$(p)$.
 
We shall present these results for Mathias
forcing only, although they can be generalized to include many more
forcing notions. One reason is that the optimal level of generality is
not clear.

\sm
 
{\bf Lemma 3.1.} {\it Let $\l P_\a ,\dot Q_\b :\a\leq \gamma
,\b<\gamma\r$ be a countable support iteration of Mathias forcing. Let
$(N,\in)$ be a countable model of ZF$^-$. Let $a^*\sub\g$ be closed such that
$a^*\in N$ and $a^*\sub N$ (so $a^*$ is countable in $V$). Let $\l
P_{a^*\cap \a}, \dot Q_\a: \a\in a^*\r$ be a countable support
iteration with domain $a^*$ of Mathias forcing. 

If $N\models p\in P_{a^*}$, there exists $q\in P_\g$ with cl$(q)=a^*$
such that $q$ is $(N,P_{a^*}, p)$-generic, that is, if $\l r_\a:\a<\g
\r$ is a sequence of $P_\g$-generic Mathias reals over $V$ with $q$
belonging
to its induced generic filter, then $\l
r_\a:\a\in a^*\r$ is $(P_{a^*})^N$-generic over $N$, with $p$ belonging to its
induced filter.}

\sm

{\it Proof:} The proof follows closely Shelah's proof [Sh{\it b}, p.90] of
preservation of properness by countable support iterations. By
induction on $j\leq \max a^*$, $j\in a^*$, we prove the following:

\sm

{\narrower{\item{($\ast$)} For every $i<j$, $i\in a^*$, for every {\bf
p} a $P_i$-name for an element of $(P_{a^*\cap j})^N\cap N$,
and for every $q\in P_i$, {\it if} $q$ is $(N,P_{a^*\cap i},$ {\bf
p}$\restrict a^*\cap i)$-generic with cl$(q)=a^*\cap i$,
{\it then} there exists $r\in P_j$ with cl$(r)=a^*\cap j$
such that $r$ is ($N,P_{a^*\cap j},$ {\bf p})-generic, and $r\restrict i=q$.

}}

\sm

\no {\it Case 1:} $j=\min a^*$. Then $P_{a^*\cap j}=\{ \emptyset \}$.
We let $r=\emptyset$.

\sm

\no {\it Case 2:} $a^*\cap j=(a^*\cap \b)\cup \{ \b\}$ for some
$\b<j$. By induction hypothesis we may assume $\b=i$. Choose $\l
r_\a:\a<i\r$ $P_i$-generic over $V$ such that $q$ belongs to the
induced generic filter. Then $\l r_\a:\a\in a^*\cap i\r$ is
$(P_{a^*\cap i})^N$-generic over $N$ with {\bf p}$[r_\a:\a <i]\restrict
a^*\cap i$ belonging to the induced filter. Hence $x:=$({\bf
p}$[r_\a:\a <i](i)) [r_\a:\a\in a^*\cap i]$ is well-defined and
$N[r_\a:\a\in a^*\cap i]\models \, ``x$ is a Mathias condition''. By
Lemma 1.2, choose a Mathias condition $y\leq x$ which is
$(N[r_\a:\a\in a^*\cap i], \name{ Q_i}[r_\a:\a\in a^*\cap i])$-generic.
In $V$ we may choose a $P_i$-name $q(i)$ for $y$ such that $q$ forces
the above to hold for $q(i)$. Then $r=q\extend \l q(i)\r$ is as
desired.

\sm

\no {\it Case 3:} $\bigcup a^*\cap j=j$. Let $\l i_n:n<\o\r$ be
increasing and cofinal in $a^*\cap j$ with $i_0=i$. Let $\l D_n:n\in
\o\r$ list all subsets of $(P_{a^*\cap j})^N$ which belong to
$N$ and are dense in the sense of $N$. We define sequences $\l
q_n:n<\o\r$ and $\l ${\bf p}$_n:n<\o\r$ such that $q_0=q$, {\bf
p}$_0=$ {\bf p}, and for all $n<\o$ the following hold:

\item{(1)} {\bf p}$_{n+1}$ is a $P_{i_n}$-name for an element of
$(P_{a^*\cap j})^N$.

\item{(2)} $q_n\in P_{i_n}$ and $q_n$ is $(N,P_{a^*\cap i_n},\hbox{
{\bf p}}_n\restrict a^*\cap i_n)$-generic.

\item{(3)} $q_{n+1}\restrict i_n=q_n$.

\item{(4)} $q_n \forces _{P_{i_n}}\, ``\hbox{{\bf p}}_{n+1}\in D_n\cap
N \hbox{ and {\bf p}}_{n+1}\leq\hbox{ {\bf p}}_n$''.

\sm

Suppose that we have already gotten $q_n$ and {\bf p}$_n$. Choose $\l
r_\a: \a<i_n\r$ $P_{i_n}$-generic over $V$ with $q_n$ belonging to its
induced generic filter. Let $s=$ {\bf p}$_n[r_\a:\a< i_n]$. Hence
$s\in (P_{a^*\cap j})^N\cap N$ by (4) in case $n>0$, and by assumption on {\bf
p}$_0$ otherwise. In $N$ we can define

$$D_n'=\{ t_0\in P_{a^*\cap i_n}:\exists t_1(t_0\extend t_1\in D_n
\hbox{ and } t_0\extend t_1\leq s)\}.$$

\no Then $N$ thinks that $D_n'$ 
is dense below $s\restrict i_n$ in $P_{a^*\cap i_n}$. 
By (2), $s\restrict i_n$ belongs to the $(P_{a^*\cap i_n})^N$-generic
filter induced by $\l r_\a:\a\in a^*\cap i_n\r$. By genericity this
filter meets $D_n'\cap N$, and hence there is $t\in D_n\cap N$ with $t\leq s$ and
$t\restrict i_n$ belonging to the filter. In $V$ we find a
$P_{i_n}$-name {\bf p}$_{n+1}$ for $t$ such that $q_n$ forces the
above properties of $t$ to hold for {\bf p}$_{n+1}$.

By induction hypothesis, $(\ast )$ is true for $i=i_n$, $j=i_{n+1}$.
Therefore there exists $q_{n+1}\in P_{i_{n+1}}$, such that (3) holds
and (2) holds for $n+1$ instead of $n$.

This finishes the construction. Now let $r=\bigcup_{n<\o}q_n$. Then
$r$ is as desired, as is easily seen.

Since $a^*$ is closed, the three cases are exhaustive. \hfill $\qed$

\sm

We start defining $(P_\g ',\leq ')$. For $\a$ an ordinal,
define $P_\a' $ as follows: 

{\narrower{\no $p\in P_\a'$ iff $p$ is a function,
dom$(p)\in [\a]^{\leq\o}$, and for all $i\in$ dom$(p)$ there exists
$u^p_i\in[i]^{\leq \o}$ such that $p(i)$ is the code of a Borel
function with domain the set of all functions $r:u^p_i\rightarrow
\fun$ and target the set of Mathias conditions.
For $i\not\in $ dom$(p)$, we let
$u^p_i=\emptyset$. 

}}

For any well-ordered set $a^*$, we can similarly define $P'_{a^*}$. If
$p\in P_\a '$, we let cl$(p)=\bigcup\{
u^p_i:i\in$ dom$(p)\}\cup$ dom$(p)$.

\sm

{\bf Remark 3.2.}  We can view $P_\g'$ as a subset of $P_\g$. Given
$p\in P_\g'$ and $i\in$ dom$(p)$, and $\l r_j:j<i\r$ $P_i$-generic
over $V$, by absoluteness we have that $p(i)\l r_j:j<u^p_i\r$ is a Mathias
condition in the extension. By the existential completeness of forcing,
there exists a $P_i$-name $\tau _i$ such that $\forces _{P_i}\, p(i)\l
\hbox{{\bf r}}_j:j\in u^p_i\r =\tau _i$. Now we can identify $p$ with  
$\l \tau_i:i<\gamma\r\in P_\g$. {\bf In the sequel we will tacitly
make use of this identification.}

\sm

We want to define a partial order $\leq '$ on $P_\g'$ such that
forcing with $(P_\g',\leq')$ will be equivalent to forcing with
$(P_\g,\leq)$. First, for $p\in P_\a'$ we define by induction on $\a\leq
\g$ when some family of reals $\l r_j:j\in u\r$ with cl$(p)
\sub u$ satisfies $p$:

\sm

\item{$\a=0$:} The only member of $P_0$ is $\emptyset$, and we
stipulate that every sequence of reals satisfies $\emptyset$;

\item{$\a=\b+1$:} $\l r_j:j\in u\r$ satisfies $p$ if $\l r_j:j\in u\r$
satisfies $p\restrict \b$ and the filter of Mathias conditions induced
by $r_\b$ contains $p(\b)\l r_j:j\in u^p_\b\r$;

\item{$\a=\bigcup\a$:} $\l r_j:j\in u\r$ satisfies $p$ if $\l r_j:j\in
u\r$ satisfies $p\restrict \b$ for all $\a<\b$.

\sm

Now let $p,q\in P_\g'$. We define:

{\narrower{\noindent $p\leq 'q$ iff dom$(q)\sub$ dom$(p)$, $u^q_i\sub
u^p_i$ for all $i\in $ dom$(p)$, and for every family of reals $\l
r_j:j\in u\r$ such that cl$(p)\sub u $ and $\l r_j:j\in u\r$
satisfies $p$, for every $i\in$ dom$(q)$ we have:

$$p(i)\l r_j:j\in u^p_i\r \leq q(i)\l r_j:j\in u^q_i\r,$$

\noindent where $\leq $ denotes the Mathias order.

}}

\sm 

Being a Borel code is a $\Pi ^1_1$ property (see [J, p. 538]).
Therefore,
by the definitions and absoluteness of $\Pi ^1_1$ statements we obtain
that the definition of $(P'_\g ,\leq ')$ is very much absolute. 

\sm

{\bf Fact 3.3.}  {\it Let $(N,\in )$ be a countable transitive model of
ZF$^-$ with $\g\in N$. Then $N\models p\in P_\g'$ iff $p\in P_\g'\cap
N$ and $N\models \hbox{ cl}(p)$ is countable. Moreover, for every
$p,q\in (P_\g')^N$ we have that $N\models p\leq 'q$ iff $p\leq 'q$.} 

\sm

Later we will use variants of this Fact without proof. In particular
we will have that $\g $ is countable in $N$. Then ``$N\models \hbox{
cl}(p)$ is countable'' follows, and we do not have to assume that
$N$ is transitive.

\sm 

We want to prove equivalence of the forcings $(P_\g,\leq)$ and
$(P_\g',\leq ')$. We start with the following easy observation:

\sm

{\bf Lemma 3.4.}  {\it If $p,q\in P_\g'$, then $p\leq' q$ implies
$p\leq q$.}

\sm

{\it Proof:} By induction on $\a\leq \g$ we prove that this is true
for $P_\a'$.

\item{$\a=0$:} clear.

\item{$\a =\b+1$:} $p\leq' q$ clearly implies $p\restrict \b\leq'
q\restrict \b$. By induction hypothesis we conclude $p\restrict \b\leq
q\restrict \b$. Let $G_\b$ be $P_\b$-generic over $V$ with $p\restrict
\b\in G_\b$. Let $\l r_j:j<\b\r$ be the sequence of Mathias reals
determined by $G_\b$. It is clear that $\l r_j:j<\b\r$ satisfies
$p\restrict \b$. By assumption we have $p(\b)\l r_j:j \in u^p_\b\r
\leq q(\b)\l r_j:j\in u^q_\b \r$. By our identification (see Remark
3.2) we have $p(\b)\l r_j:j \in u^p_\b\r =p(\b)[G_\b]$ and
$q(\b)\l r_j:j\in u^q_\b \r =q(\b)[G_\b]$. Consequently $p\restrict \b
\forces _{P_\b}\, p(\b)\leq q(\b)$, and hence $p\leq q$.

\item{$\a=\bigcup\a$:} clear by induction hypothesis and definition of
the partial orders. \hfill $\qed$

\sm

The next lemma shows that $P_\g '$ is a dense subset of $P_\g$. In the
proof  we will use the following coding of Mathias conditions by reals
$x\in \fun$ with the property $\forall i,j(0<i<j\Rightarrow
x(i)< x(j))$: such $x$ codes the Mathias condition
$(\hbox{ran}\, x\restrict [1,x(0)),\hbox{ ran}\, x\restrict [x(0),
\infty ))$. Hence we may assume that a $P_i$-name for a Mathias
condition is a sequence $\l f_n: n<\o\r$ such that $f_n:
A_n\rightarrow \o$, where $A_n$ is a countable antichain of $P_i$.

For $p\in P_\g$ and sequence of reals $\bar r=\l r_j:j\in u\r$ with
cl$(p)\sub u$, we define by induction on $i\leq \g$, $i\in\hbox{
dom}(p)$, 

\item{(a)} $\bar r$ {\it evaluates} $p(i)$;

\item{(b)} $p(i)[\bar r]$, if $\bar r$ evaluates $p$.

\sm

\no {\it Case 1:} $i=0$. $\bar r$ evaluates $p(i)$, $p(i)[\bar r]=p(i)$.

\sm

\no {\it Case 2:} $i>0$. Then $p(i)=\l f_n:n<\o\r$, where
$f_n:A_n\rightarrow \o$ and $A_n\sub P_i$ is a countable antichain.
We define that $\bar r$ evaluates $\g$ if:

\item{(1)} for every $n<\o$, every $q\in A_n$, and every $\b\in\hbox{
dom}(q)$, $\bar r$ evaluates $q(\b)$;

\item{(2)} for every $n<\o$ there exists a unique $q\in A_n$ such that
for all $\b\in\hbox{ dom}(q)$, $q(\b)[\bar r]$ belongs to the filter
on $Q$ induced by $r_\b$;

\item{(3)} the real $x$ defined by $x(n)=f_n(q)$, where $q\in A_n$ is
the unique member as in (2), codes a Mathias condition (i.e. $\forall
i,j( 0<i<j<\o \Rightarrow x(i)<x(j))$).

If (1)--(3) hold, $p(i)[\bar r]$ is defined as the Mathias condition
coded by $x$.

\sm

The set of sequences $\bar r=\l r_j: j\in\hbox{ cl}(p(i))\r$ which
evaluate $p(i)$ is a Borel set with code $p(i)$; for it is not
difficult, though tedious, to show that it has a $\Delta
^1_1(p(i))$-definition (see [JSp], where the details are worked out).
First, $\bar r$ evaluates $p(i)$ iff there exists a sequence of reals
which are the evaluations by $\bar r$ of all the names which belong to
the transitive closure of $p(i)$, such that $p(i)$ can be evaluated
from these using $\bar r$. Since $p(i)$ is hereditarily countable
there is only one existential real quantifier, and the others are
number quantifiers. Second, if such a sequence of reals exists, then
it is unique, hence we can turn this statement into a universal
statement. Now by Suslin's Theorem (see [J, p.502]) we are done.
 
By a similar argument, the map sending $\bar r$, which evaluates $p(i)$,
to $p(i)[\bar r]$ has a Borel definition.

\sm

{\bf Lemma 3.5.}  {\it For every $p\in P_\g$ there exists $p'\in P_\g '$
such that $p'\leq p$.}

\sm

{\it Proof:}  For each $i\in$ dom$(p)$ let $u^i_p= $ cl$(p(i))$. Then
$u^i_p$ is countable. We define $p'(i): \{ \bar r: \bar r: u^p_i
\rightarrow \fun \} \rightarrow Q$ ($Q$ is Mathias forcing) by cases
as follows: If $\bar r$ evaluates $p(i)$, we let $p'(i)(\bar
r)=p(i)[\bar r]$, otherwise we let $p(i)(\bar r)$ be the maximum
element of $Q$. By the remarks above, $p'(i)$ is a total Borel
function as desired. Now let $p'=\l p'(i):i\in\hbox{ dom}(p)\r$. Then
clearly $p'\in P'_\g$. By induction on $i\in \hbox{ dom}(p')$ it is
easy to prove that if $\bar r=\l r_j:j<i\r$ is $P_i$-generic over $V$
and contains $p'\restrict i$ in its generic filter, then $\bar r$
evaluates $p(i)$ and $p'(i)(\bar r)=p(i)[\bar r]$; hence $p'\restrict
i\forces _{P_i}\, p'(i) =p(i)$.\hfill $\qed$

\sm

In order to conclude that forcings $(P_\g ,\leq)$ and $(P_\g ',\leq
')$ are equivalent it is enough to prove the following:

\sm

{\bf Lemma 3.6.} {\it For all $p,q\in P'_\g$ with $p\leq q$ there exists
$r\in P'_\g$ with $r\leq 'p$ and $r\leq 'q$.}

\sm

{\bf Corollary 3.7.}  {\it Forcings $(P_\g ,\leq)$ and $(P_\g ',\leq
')$ are equivalent.} 

\sm

{\it Proof of 3.7:} By Lemma 3.5 it is enough to show that $(P'_\g ,\leq)$ and
$(P'_\g ,\leq ')$ are equivalent. Let $D$ be dense open in
$(P'_\g,\leq)$, and let $p\in P'_\g$. Let $q\in D$, $q\leq p$. By
Lemma 3.6 there is $r\in P'_\g$ with $r\leq 'p$ and $r\leq 'q$. By 3.4
we have $r\leq q$, and hence $r\in D$. Therefore $D$ is dense in
$(P'_\g, \leq ')$. Conversely, if $D$ is dense in $(P'_\g,\leq ')$, 
then $D$ is dense in $(P'_\g ,\leq )$ by Lemma 3.4.

\relax From Lemma 3.6 it follows that for all $p,q\in P'_\g$, $p,q$
are incompatible with respect to $\leq $ iff they are incompatible
with respect to $\leq '$. Therefore every $(P'_\g, \leq)$-name is a
$(P'_\g, \leq ')$-name and vice versa.

It follows that if $G$ is a $(P'_\g ,\leq)$-generic filter, then $G$
is also $(P'_\g ,\leq ')$-generic, and if $G'$ is $(P_\g ',\leq
')$-generic, then $G=\{ p\in P_\g ':\exists q\in G'(q\leq p)\}$ is
$(P_\g ',\leq ')$-generic, and then $V[G]=V[G']$.\hfill $\qed$

\sm

The following will be crucial for proving Lemma 3.6:

\sm

{\bf Lemma 3.8.} {\it Let $a^*$ be a countable closed set of ordinals,
and let $p\in P_{a^*}'$. Let $(N,\in )$ be a countable elementary
substructure of $(H(\chi) ,\in )$ for some large enough regular $\chi
$, such that $p, a^*\in N$. There exists $q\in P_{a^*}'$, $q\leq 'p$,
such that for every sequence of reals $\bar r=\l r_l:l\in a^*\r$ which
satisfies $q$, $\bar r$ is $(P_{a^*},\leq )$-generic over $N$.}

\sm

{\it Proof:} By induction on $j\in a^*$ we prove the following:

\sm 

{\narrower{\item{$(\ast )$} For every $i<j$, $i\in a^*$, for every
$P_{a^*\cap i}$-name {\bf p} for a member of $N\cap P_{a^*\cap j}$, and for
every $q\in P'_{a^*\cap i}$, {\it if} every sequence of reals $\bar
r=\l r_l:l\in a^*\cap i\r$ which satisfies $q$ is $P_{a^*\cap
i}$-generic over $N$, and $q\forces _{P_{a^*\cap i}}\, \hbox{{\bf
p}}\restrict i\in G_{a^*\cap i}$, {\it then} there exists $r\in
P_{a^*\cap j}'$ such that $r\restrict a^*\cap i=q$, every $\l r_l:l\in
a^*\cap i\r$ which satisfies $r$ is $P_{a^*\cap j}$-generic over $N$,
and $r\forces _{P_{a^*\cap j}}\, \hbox{{\bf p}}\in
G_{a^*\cap j}$.

}}

\sm

\no {\it Case 1:} $j=\min a^*$. Let $r=\emptyset$.

\sm

\no {\it Case 2:} $a^*\cap j=(a^*\cap \b)\cup \{ \b\}$ for some
$\b<j$. By induction hypothesis we may assume $\b=i$. Let $\bar r=\l
r_l:l\in a^*\cap i\r$ satisfy $q$. By assumption, $\bar r$ is
$P_{a^*\cap i}$-generic over $N$ and {\bf p}$[\bar r]\restrict a^*\cap
i$ belongs to the generic filter induced by $\bar r$. By absoluteness,
$x:=$ ({\bf p}$[\bar r](i))[\bar r]$ is a Mathias condition in $V$,
say $x=(u^x,a^x)$. Using $N[\bar r]$ as a code we may effectively
construct $u\in \inf$ which is Mathias-generic over $N[\bar r]$ with
$x$ belonging to the generic filter induced by $u$. Let $y=(u^x,
a^x\cap u)$. Then every real $r_i$ which satisfies $y$ is
Mathias-generic over $N[\bar r]$ (see Lemma 1.2). Moreover, the
function sending $\bar r$ to $y$ is
Borel. Denote it with $r(i)$. Then we may let $r=q\extend \l r(i) \r$.

\sm
\no {\it Case 3:}  $a^*\cap j$ is unbounded in $N\cap j$. We choose 
$\l i_n:n<\o\r$ increasing and cofinal in $N\cap j$ with $i_0=i$. Let $\l
D_n:n<\o\r$ list all dense subsets of $P_{a^*\cap j}$ in $N$. We define two
sequences $\l q_n:n<\o\r$ and $\l${\bf p}$_n:n<\o\r$ such that
$q_0=q$, {\bf p} $=$ {\bf p}$_0$, and for all $n<\o$ the following
hold:

\sm

\item{(1)} {\bf p}$_{n+1}$ is a $P_{i_n}$-name for a member of
$P_{a^*\cap j}\cap N$;

\sm

\item{(2)} $q_n\in P'_{a^*\cap i_n}$, and for every $\bar r=\l r_l:l\in
a^*\cap i_n\r$ which satisfies $q_n$, $\bar r$ is $P_{a^*\cap
i_n}$-generic over $N$, and $q_n\forces _{P_{a^*\cap i_n}}\,
\hbox{{\bf p}}_n\restrict a^*\cap i_n\in G_{a^*\cap i_n}$;

\sm

\item{(3)} $q_{n+1}\restrict i_n=q_n$;

\sm

\item{(4)} $q_n\forces _{P_{a^*\cap i_n}}\, \hbox{{\bf p}}_{n+1} \in
D_n\cap N \hbox{ and } \hbox{{\bf p}}_{n+1}\leq \hbox{ {\bf p}}_n $.

\sm

The construction is analogous to the proof of Lemma 3.1. 

Now let $r=\bigcup_{n<\o}q_n$, and let $\bar r=\l r_l:l\in a^*\cap
j\r$ satisfy $r$. We have to show that $\bar r$ is $P_{a^*\cap
j}$-generic over $N$. Let $G\sub P_{a^*\cap j}$ be the filter induced
by $\bar r$. Then $r\in G$. We have to show that $D_n\cap G\ne
\emptyset$ for all $n<\o$.  Let $n<\o$. We claim that $p_{n+1}:=$ {\bf
p}$_{n+1}[\bar r\restrict i_n]\in G\cap D_n$. By (2) and (3), $\bar
r\restrict i_n$ is $P_{a^*\cap i_n}$-generic over $N$, and hence
$p_{n+1}\in D_n$ by (4).  To prove $p_{n+1}\in G$ it is enough to show
that $p_{n+1}\restrict i_m\in G_{a^*\cap i_m}$ for all
$n<m<\o$. For this, by induction on $m$ show (using (4)) that $p_m\leq
p_{n+1}$. This suffices, since by (2), $p_m\restrict a^*\cap i_m\in
G_{a^*\cap i_m}$. This finishes the proof of $(\ast )$.

Applying $(\ast )$ for $i=\min (a^*)$ and $j=\max (a^*)$, we get $q\in
P'_{a^*}$ such that every $\bar r=\l r_l:l\in a^*\r$ which satisfies
$q$ is $(P_{a^*},\leq )$-generic over $N$ and contains $p$ in its
induced filter. We have to show that $q\leq 'p$. By contradiction,
suppose $\bar r=\l r_l:l\in a^*\r$ satisfies $q$ and there is $i\in $
dom$(q)$ such that $q(i)\l r_l:l\in a^*\cap i\r \not\leq p(i)\l
r_l:l\in u^p_i\r$. We can choose $r_i'$ which satisfies $q(i)\l
r_l:l\in a^*\cap i\r $, but not $p(i)\l
r_l:l\in u^p_i\r$. Choose $\l r_l':l\in a^*\setminus (i+1)\r$
arbitrary such that $\bar r':=\l r_l:l\in a^*\cap i\r\extend \l
r_l':l\in a^*\setminus i\r$ satisfies $q$. By the above, $\bar r'$ is
$P_{a^*}$-generic over $N$, containing $p$ in its generic filter. But
this is impossible by the choice of $r_i'$. \hfill $\qed$

\sm

We are now able to give the proof of Lemma 3.6.

\sm

{\it Proof of 3.6:} Let $p,q\in P'_\g$ with $P_\g \models p\leq q$.
Let $a^*=$ cl$(p)$. Hence we have $p,q\in P'_{a^*}\sub P_{a^*}$. 
We need the following claim:

\sm

{\bf Claim:} $P_{a^*}\models p\leq q$.

\sm

{\it Proof of the Claim:} Otherwise, let $i\in $ dom$(p)$ be minimal
such that $\lnot (p\restrict i \forces _{P_{a^*\cap i}}$ $p(i)$ $\leq
q(i))$. Choose $r\in P_{a^*\cap i}$ such that $P_{a^*\cap i}\models
r\leq p\restrict i$ and $r\forces _{P_{a^*\cap i}}\, p(i)\not\leq
q(i)$.

Let $(N,\in )$ be a countable elementary substructure of $(H(\chi),\in
)$, $\chi$ large enough and regular, containing everything relevant.
By Lemma 3.1 there exists $q_1\in P_i$ which is $(N, P_{a^*\cap i},
r)$-generic. Let $\bar r =\l r_j:j<i\r$ be $P_i$-generic
over $V$ with $q_1$ belonging to the induced filter.  Then $\l
r_j:j\in a^*\cap i\r$ is $P_{a^*\cap i}$-generic over $N$, with
$r$ belonging to the induced filter. We conclude that on
the one hand, $V[r_j:j<i]\models p(i)[r_j:j<i]\leq q(i)[r_j:j<i]$, but
on the other hand, $N[r_j:j\in a^*\cap i]\models p(i)[r_j:j\in a^*\cap
i]\not\leq q(i)[r_j:j\in a^*\cap i]$. But $p(i)[r_j:j<i]=p(i)[r_j:j\in
a^*\cap i]$, and similarly for $q(i)$. Since the Mathias order is
absolute, we have a contradiction. \hfill $\qed$

\sm

Let $(N,\in )$ be as in the proof of the Claim. By Lemma 3.8, there
exists $r\in P'_{a^*}$ with $r\leq 'p$ such that every sequence of
reals $\bar r=\l r_j:j\in a^*\r$ which satisfies $r$ is
$P_{a^*}$-generic over $N$. Given such $\bar r$ and $i\in $ dom$(p)$,
$p\restrict i$ belongs to the generic filter on $P_{a^*}\cap N$
induced by $\bar r\restrict a^*\cap i$, and hence by the Claim, $N[\bar
r\restrict a^*\cap i] \models p(i)[\bar r\restrict a^*\cap 
i]\leq q(i)[\bar r\restrict
a^*\cap i]$. But $p(i)[\bar r\restrict a^*\cap
i]=p(i)(r_j:j\in u^p_i)$, and similarly
for $q$. By absoluteness of the Mathias order and by $r\leq 'p$ we
obtain $r(i)(r_j:j\in u^r_i)\leq p(i)(r_j:j\in u^p_i)\leq
q(i)(r_j:j\in u^q_i)$. Since $\bar r$ and $i$ were arbitrary we
conclude that $r\leq 'q$. \hfill $\qed$
   
\sm

The proof of Corollary 3.7 now being complete, {\bf throughout the
rest of this paper we identify $(P_\g ,\leq)$ with $(P_\g ',\leq ')$.}

\sm

{\bf Definition 3.9.} If $u\sub \g$ is finite and $p,q\in P_\g$, then $q\leq
_up$ is defined by: $q\leq p$ and for all $\a\in u$,
$q\restrict \a\forces _{P_\a}\, ``q(\a)\hbox{ and }p(\a)\hbox{ have
the same first coordinate''}$.

\sm

By arguments which are standard by now, we obtain the following Lemma.
Note that it makes sense only in the light of Corollary 3.7. For the
proof, make a similar inductive construction as we did now several
times. At successor steps use Lemma 1.2 to get generic conditions
which are pure extensions, if required by $u$.

\sm

{\bf Lemma 3.10.} {\it Let $(N,\in)$ be a countable model of ZF$^-$
such that $\g$ is countable in $N$. If $p\in P_\g\cap N$, and $u\in
[\g ]^{<\o}$, there exists $q\in P_\g$ such that $q\leq _u p$ and $q$
is $(N,P_\g)$-generic.}

\sm

For the proof that potential counterexamples to Propositions 2.3 and
2.4 are added by an iteration of countable length, we will also need the
following lemma.

\sm

{\bf Lemma 3.11.} {\it Suppose $a^*\sub \g$ is a countable closed set
of ordinals, $P_{a^*}$ is a countable support iteration of Mathias
forcing with domain $a^*$, and $p\in P_{a^*}$.  Let $(N,\in )$ be a
countable model of ZF$^-$ with $\g\in N$, and suppose that $a^*\sub
N$, $a^*\in N$, $p\in N$, and $N\models p\in P_{a^*}$.

There
exists $q\in P_{a^*}$ and a $P_{a^*}$-name ${\bar{\hbox{{\bf r}}}}'_\g
=\l{\hbox{{\bf r}}}'_l:l<\g\r $ such that $q\leq p$ and, letting
${\bar{\hbox{{\bf r}}}}_{a^*}=\l\hbox{{\bf r}}_l:l\in a^*\r $ be a name
for the $P_{a^*}$-generic sequence of Mathias reals, we have

$$q\forces_{P_{a^*}}\, ``\bar {\hbox{{\bf r}}}'_\g \hbox{ is
}P_\g\hbox{-generic over }N,\hbox{ and }\forall l\in a^*(\hbox{{\bf
r}}'_l =
\hbox{{ \bf r}}_l)\hbox{''}.$$}

\sm

{\it Proof:} By induction on $j\leq \g$, $j\in N$, we prove the following:

{\narrower{\item{($\ast $)} Suppose $i\in j$, $i\in N$, $q\in P_{a^*\cap i}$,
and ${\bar{\hbox{{\bf r}}}}'_i =\l\hbox{{\bf r}}'_l:l<i\r $ is a
$P_{a^*\cap i}$-name such that $q\leq p\restrict a^*\cap i$ and 

$$q\forces_{P_{a^*\cap i}}\, {\bar {{\hbox{{\bf r}}}}}'_i \hbox{ is
}P_i\hbox{-generic over }N \hbox{ and }\forall l\in a^*\cap
i({\hbox{{\bf r}}}'_l={\hbox{ {\bf r}}}_l).$$

Then there exists $r\in P_{a^*\cap j}$ and ${\bar{\hbox{{\bf
r}}}}'_j =\l\hbox{{\bf r}}'_l:l<j\r $ such that $r\restrict a^*\cap
j=q$, 

$r\leq p\restrict a^*\cap j$,
${\bar{\hbox{{\bf r}}}}'_j \restrict i={\bar{\hbox{{\bf
r}}}}'_i$, and

$$r\forces_{P_{a^*\cap j}}\, ``{\bar{\hbox{{\bf r}}}}'_j \hbox{ is
}P_j\hbox{-generic over }N\hbox{ and }\forall l\in a^*\cap j(
\hbox{{\bf r}}'_l=\hbox{{\bf r}}_l) .$$

}}

\sm

\no {\it Case A:} $N\cap j=(N\cap \b)\cup \{\b\}$, for some $\b<j$: Then $j=\b +1$, since
$N\models $ ZF$^-$, and so $\b+1\in N$. Hence we may assume $\b=i$.
\sm

\no {\it Case A1:} $i\in a^*$. Let $\bar r_{a^*\cap i}=\l r_l:l\in
a^*\cap i\r$ be $P_{a^*\cap i}$-generic over $V$ with $q$ in its
generic filter. Let $\bar r'_i={\bar{\hbox{{\bf r}}}}'_i[\bar r_{a^*\cap
i}]$. Then $N[\bar r'_i]\in V[\bar r_{a^*\cap i}]$ and $N[\bar
r_i']\models \hbox{ ZF}^-$. By assumption we have $p(i)[\bar r_{a^*\cap
i}]= p(i)[\bar r'_i]$. Let $x$ be this common value. Then $x$ is a
Mathias condition. By Lemma 1.2, in $V[\bar r_{a^*\cap i}]$ we may choose
a Mathias condition
$y\leq x$ such that every $z\in \inf$ with $u^y\sub z\sub u^y\cup a^y$
is Mathias generic over $N[\bar r'_i]$. In $V$ we have a $P_{a^*\cap
i}$-name $q_i$ such that $q$ forces that all the above holds for $q_i$
instead of $y$. Now let $r=q\extend \l q_i\r$ and $\hbox{{\bf r}}'_i=
\hbox{ {\bf r}}_i$.

\sm

\no {\it Case A2:} $i\not\in a^*$. Then $P_{a^*\cap j}=P_{a^*\cap i}$.
Since $N$ is countable, in $V$ there exists a $P_{a^*\cap i}$-name
$\hbox{{\bf r}}_i'$ such that $q$ forces that $\hbox{{\bf r}}_i'$ is
Mathias generic over $N[{\bar{\hbox{{\bf r}}}}'_i]$. We let $r=q$ and
${\bar{ \hbox{{\bf r}}}}_j'={\bar{\hbox{{\bf r}}}}_i'\extend
\l\hbox{{\bf r}}_i'\r$.

\sm

\no {\it Case B:} $N\cap j$ is unbounded in $N\cap j$:
\sm

\no {\it Case B1:} $j\in a^*$. Since $a^*$ is closed and $a^*\sub N$,
we conclude that either $a^*\cap j$ is bounded in $a^*\cap j$, or else
$a^*\cap j$ is unbounded in $j$. In the first case we may assume
$i>\max (a^*\cap j)$, and proceed as in Case A2. In the latter
case, a similar diagonalization as in 3.1 and 3.8 works.

\sm

\no {\it Case B2:} $j\not\in a^*$. Since $a^*$ is closed, $a^*\cap j$ is
bounded below $j$. Hence we may assume $i>\max (a^*\cap j)$. Then
$P_{a^*\cap j}=P_{a^*\cap i}$, and as in Case A2, in $V$
there exists a $P_{a^*\cap i}$-name $\l \hbox{{\bf r}}_l':i\leq l<j\r$
such that $q$ forces that  $\l \hbox{{\bf r}}_l':i\leq l<j\r$ is
$P_j$/${\bar {\hbox{{\bf r}}}}_i'$-generic over $N$. We let $r=q$ and $
{\bar {\hbox{{\bf r}}}}_j'= {\bar {\hbox{{\bf r}}}}_i'\extend \l \hbox{{\bf
r}}_l':i\leq l<j\r$. \hfill $\qed$

\bigskip

{\hgross 4. Proof of Proposition 2.3}

\bigskip

The following Lemma will give us the $\o _1$-club for Proposition 2.3.

\sm

{\bf Lemma 4.1.}  {\it Suppose $V\models CH$. Let $\iter $ be a
countable support
iteration of Mathias forcing. Let $G_{\o _2}$ be $P_{\o _2}$-generic over $V$
and, for $\d <\o _2$, $r_\d $ the $\name{ Q_\d} [G_\d ]-$generic 
real determined by
$G_{\o _2}$. Then the set $S$ of $\d \in S^2_1$ such 
that for some $\a _\d < \d $

$$\P ^{V[\{ G_{\a _\d }, r_\d \}]} = \P ^{V[G_{\d +1}]}\eqno(\ast )
$$ 

\noindent is nonstationary.}

\sm

{\it Proof:}  Suppose that $S$ is stationary. We will derive a contradiction.
For $\d \in S$ choose $p_\d \in P_{\d +1}$ forcing ($\ast $). Since $\d \in
S^2_1$ and $p_\d $ is hereditarily countable, without loss
of generality we may assume that $p_\d (\d )$ is a $P_{\a _\d} $-name and
$\sup (\dom (p_\d \restrict \d )) <\a _\d $. Otherwise increase $\a _\d $, and then 
$(\ast )$ still holds of course. By Fodor's Theorem and $V[G_\a ]\models CH$
for $\a <\o _2$, there exist $\a ^* < \o _2 $, $p\in P_{\a ^*}$ and a 
stationary $S_1\sub S$ such that $\forall \d \in S_1 (\a _\d =\a ^* \wedge
p_\d \restrict \d =p)$. Hence in $V[G_{\a ^*}]$ we can compute $p_\d (\d )[G_\d ]$
for $\d \in S_1$. Again by the CH in $V[G_{\a ^*}]$ and the $\aleph
_2$-completeness of the nonstationary ideal on $\o _2$,
there exist a stationary $S_2\sub S_1$ and $q\in
\name{ Q_{\a ^*}}[G_{\a ^*}]$ such 
that $\forall \d \in S_2 (p_\d (\d )[G_\d ]=q)$.

Let $G(\o _2)$ be $Q ^{V[G_{\o _2}]}$-generic over $V[G_{\o
_2}]$, where $Q$ is Mathias forcing,  
such that $q\in G(\o _2)$. Let $r_{\o _2}$ be the corresponding Mathias
real, and let $G_{\o _2 +1} =G\ast G(\o _2)$. By Theorem 2.1, ${\cal
P}(\o)$/fin is $\aleph _1$-distributive in $V[G_{\o _2}]$. Since
Mathias forcing is the composition of ${\cal
P}(\o)$/fin and some $\sigma$-centered forcing, it follows that 
$V[G_{\o _2
+1}]\models {\frak c}=\o _2$. By properness and $V\models CH$ we have
$V[G_{\a ^*}, r_{\o _2}]\models CH$. (If you do not see this, 
let $V=L$ and use [J, 15.3., p.130].) Hence there exists $\a ^* <\a  <
\o _2$ such that $r_\a
\not\in V[G_{\a ^*}, r_{\o _2}]$. Hence in $V[G_{\o _2}]$ there exists $q_1\in 
Q ^{V[G_{\o _2}]} \cap G(\o _2)$, $q_1\leq q$, forcing this. Let 
$\a < \gamma <\o _2$ such that $q_1\in V[G_\gamma ]$. 
By genericity there exists $\d \in S_2\cap [\gamma , \o _2) $ such that,
if $q_1=( u,a) $ then $u\sub r_\d \sub u\cup a$, that is, $q_1$ belongs
to the generic filter generated by $r_\d $. Let $q_2=( u, a\cap r_\d ) $.
Then $q_2\in Q ^{V[G_{\o _2}]}$ and $q_2\leq q_1$. 

Let $r$ be $Q ^{V[G_{\o _2}]}$-generic over $V[G_{\o _2}]$ such
that $u\sub r\sub u\cup (a\cap r_\d )$. Then $r$ is an infinite subset
of $r_\d $. By the remark preceding Lemma 1.2, we have that $r$ is $Q
^{V[G_\d ]}$-generic over $V[G_\d ]$.  From $(\ast )$ and the choice
of $q$ we conclude that $r_\a \in V[G_{\a ^*},r]$. But on the other
hand, $q_1 $ belongs to the generic filter induced by $r$, and we
conclude $r_\a \not\in V[G_{\a ^*}, r]$, a contradiction. \hfill
$\square $

\sm

Let $C\sub S^2_1\setminus S$ be $\o _1$-club, where $S$ is as
in Lemma 4.1. We claim that $C$ serves for Proposition 2.3. By
contradiction, suppose that this is false. Hence there exist
$\a \in C$, $p^*\in P_{\o _2 }/\name{G_\a} $, and $\name{ r}$ such that 

$$p^* \forces _{P_{\o _2 }/\name{G_\a}}\, \name{ r} \hbox{ induces a
Ramsey ultrafilter on }(\inf)^{V[\name{G_\a}]}\hbox{ and }\name{
r}\not\in V[\Name{G_{\a +1}}].\eqno{(\ast )}$$

Since forcing $P_{\o _2}/\name{G_\a} $ is equivalent to a countable
support iteration of length $\o _2 $ of Mathias forcing in $V[G_\a ]$
(see [B, $\S 5$]),
for notational simplicity we assume $\a =0$ for the moment, and later we
will remember that really $V=V[G_\a ]$ for some $\a \in C$ and derive a
final contradiction. 

First we show that by the absoluteness results from $\S 3$ we may
assume that $\dot r$ is added by an iteration of countable length. Let
$a^*=$ cl$(p^*)$. So $a^*\sub \o_2$ is countable. We may assume that
$0\in a^*$ and $a^*$ is closed.

\sm

{\bf Lemma 4.2} {\it Assuming $(\ast )$, it is true that} $p^* \forces _{P_{a^* }}\, \name{r}$ {\it induces a Ramsey ultrafilter on} $(\inf )^{V}$ {\it and} $
\name{ r} \not\in V[G_0]$.
 
\sm

{\it Proof:} (a) $p^*\forces _{P_{a^* }}\, \name{ r} \hbox{ induces an
ultrafilter on }(\inf)^{V}$: Otherwise there exists $a\in (\inf)^V$
and $p\in P_{a^*}$ such that $p\leq p^*$ and $p\forces _{P_{a^*}}\,
``\name{ r}\cap a\hbox{ and }\name{ r}\cap (\o\setminus a)$ are both
infinite.'' Let $\chi $ be large enough and regular, and let
$(N,\in)\prec (H(\chi ),\in)$ be countable, containing everything
relevant. By Lemma 3.1 choose $q\in P_{\o _2}$ such that $q$ is
$(N,P_{a^*}, p)$-generic, and let $\l r_\a :\a\in \o _2\r$ be $P_{\o
_2}$-generic over $V$, with induced filter $G$, such that $q\in G$.
Then $\l r_\a :\a\in a^*\r$ is $P_{a^*}$-generic over $N$ with $p$,
and hence also $p^*$, in its generic filter, denoted $G_{a^*}$. Then
clearly $p^*\in G$. We obtain that $V[G]\models ``\name{ r}[G]\sub ^* a$
or $\name{ r}[G]\sub ^*\o\setminus a$'', and $N[G_{a^*}]\models |\name{
r}[G_{a^*}] \cap a |=|\name{ r}[G_{a^*}]\cap (\o\setminus a)|=\o$. But
clearly $\name{r}[G]=\name{ r}[G_{a^*}]$, a contradiction.

(b) $p^*\forces _{P_{a^*}}\, \name{ r}\not\in V[G_0]$: Otherwise there is
$p\in P_{a^*}$, $p\leq p^*$, such that $p\forces _{P_{a^*}}\, \name{
r}\in V[G_0]$. Choose $(N,\in)$, $q$ and $G$ as in (a), and let
$G_{a^*}$ be defined as there. Then $\name{ r}[G]\not\in V[G_0]$, $\name{
r}[G_{a^*}]\in N[G_0]$ and $\name{ r}[G]=\name{ r}[G_{a^*}]$. Since
$N[G_0]\in V[G_0]$, we have a contradiction.

(c) $p^*\forces _{P_{a^*}}\, \name{ r} $ induces a {\it Ramsey}
ultrafilter on $(\inf )^{V}$: Otherwise there exist $p\in P_{a^*}$ and
$f\in (\fun )^V$ such that if $\D$ is a $P_{a^*}$-name for the filter
induced by $\name{ r}$ we have that $p\forces _{P_{a^*}}\, f$ is
unbounded but not one-to-one modulo $\D$. Let $(N,\in )$ be as above
containing everything relevant. We can get $q\in P_{a^*}$, $q\leq p$,
as in Lemma 3.11. Let $\bar r_{a^*}=\l r_l:l\in a^*\r$ be
$P_{a^*}$-generic over $V$ containing $q$ in its generic filter. By
Lemma 3.11, in $V[\bar r_{a^*}]$ there exists $\bar r_{\o _2}'=\l
r_l':l<\o _2\r$ such that $\bar r_{\o_2}'$ is $P_{\o_2}$-generic over
$N$ and $r_l=r_l'$, for all $l\in a^*$. We obtain that $\name{ r}[\bar
r_{\o_2}']= \name{ r}[\bar r_{a^*}]$. Let $r$ be the common value. Then
$r$ induces the same filter, say $D$, in $V[\bar r_{a^*}]$ and in
$N[\bar r_{\o _2}']$, and also $f$ is unbounded modulo $D$ in both
models. Hence by construction, on the one hand we have that $V[\bar
r_{a^*}] \models f$ is not one-to-one modulo $D$, but one the other hand
$N[\bar r_{\o _2}']\models f$ is one-to-one modulo $D$. Since $N[\bar
r_{\o _2}']\in V[\bar r_{a^*}]$ we have a contradiction. \hfill $\qed$

\sm

Continuing the proof of Proposition 1, let $\delta =$ o.t.$(a^*)$.
Then $\delta <\omega _1$,
and clearly $P_{a^*}$ and $P_\delta $ are isomorphic. Then our
assumption $(\ast )$ becomes:

$$p^* \forces _{P_{\delta }}\, \name{ r} \hbox{ induces a Ramsey
ultrafilter on }(\inf)^{V}\hbox{ and }\name{ r}\not\in V[G_0].\eqno{(\ast
\ast )}$$

Let $\D $ be a $P_\delta$-name for the filter on $(\inf )^V$ induced
by $\name{ r}$. 
In $V$, let $( N,\in ) $ be a countable elementary substructure of
$( H(\chi ),\in ) $, where $\chi $ is a large enough regular cardinal,
such that $\delta ,p^* , \D , \name{ r}
\in N$. {\bf This $N$ will be fixed for the rest of this section.} Let
$G_0$ be $Q_0$-generic, containing a $( N,Q_0) $-generic
condition below $p^*(0)$. In $V[G_0]$ we define:

$$\Y =\{ Y:\exists ( N[G_0], P_{\delta}/\name{G_0} ) \hbox{-generic }q(q\leq
p^*\restrict [1,\d )   \wedge q \Vdash _{P_\d/\name{G_0}} ``\D \cap N =Y\hbox{''})\} .$$

Since every Ramsey ultrafilter is a $p$-point (see $\S 1$), and
every $Y\in \Y $ is a countable subset of the denotation of $\D $
in a $P_{\delta }/\name{G_0}$-generic extension of 
$V[G_0]$, and $\D $ is forced to be
a Ramsey ultrafilter on $(\inf )^V$, we conclude that such $Y$ is
definable from $([\o ]^\o )^N$ and a member of $([\o ]^\o )^V $, and
hence $\Y \sub V$.


\sm

{\bf Lemma 4.3.} {\it $\Y $ is a $\Sigma ^1_2$ set in $V[G_0]$.}

\sm

{\it Proof:}  We show that $Y\in \Y $ is equivalent to saying:

\sm

{\narrower{\no There exists a countable model $( M,\in )$ such that
$N[G_0]\cup \{ N[G_0], Y\} \sub M$, $( M,\in )$ $ \models ZF^-$, and
$( M,\in ) \models \exists q\in P_{\delta }/\name{G_0}(q$ is $(
N[G_0], P_{\delta }/\name{G_0} )$-generic and $ q \Vdash
_{P_\d/\name{G_0}} ``\D \cap (\inf ) ^N =Y\hbox{''}) .$

}}

\sm

It is well-known (see [J, the proof of 41.1., pp.527f.]) that the
quantification over countable models as above is equivalent to
quantifying over structures $( \o ,R) $, where $R$ is a well-founded
binary relation, -- which makes the formula no worse (and no better)
than $\Sigma ^1_2
$ --, and that the rest is arithmetical.  

If $Y\in \Y $, then choosing a countable $( M,\in ) $ which is
elementarily embeddable into 
$( H(\chi )^{V[G_0]},\in ) $ and contains $N[G_0]\cup \{ N[G_0],Y\} $,
we easily see that one implication holds.

Conversely, if $( M,\in ) ,Y,q $ are given as above, then by Lemma
3.10, in $V[G_0]$ choose $q_1\leq q$ which is $( M,P_{\delta }/\name{G_0} )
$-generic. Here we use again the fact that $P_\delta /\name{G_0}$ is
equivalent to a countable support iteration of Mathias forcing.  Then
clearly $q_1 $ is also $( N[G_0], P_{\delta }/\name{G_0}) $-generic, and
$q_1
\Vdash _{P_\delta /\name{G_0} } ``\D
\cap (\inf )^N=Y $'' holds in $V[G_0]$. In fact, let
$G_{1}$ be $P_\delta /\name{G_0}$-generic over $V[G_0]$, containing $q_1$.
Then $G_{1}$ is $P_\delta /\name{G_0}$-generic over $M$ and contains $q$.
By assumption on $M$, $G_{1}$ is $P_\delta /\name{G_0}$-generic over
$N$.
Moreover, $\name{r}[G_0\ast G_{1}]$ is the same real in $V[G_0\ast
G_{1}]$ and $N[G_0\ast G_{1}]$. Hence we are done.
\hfill $\square $

\sm

The crucial fact, whose proof will require considerable space, is that
${\cal Y} $ is uncountable. Then we obtain that in
$V[G_0]$, ${\cal Y}$ is an uncountable $\Sigma ^1_2$ set which is a subset
of $V$. By a well-known result of descriptive set theory (see the
remark after Corollary 4.10, below), either
${\cal Y}$ has a perfect subset, or else ${\cal Y}$ is the union of
$\aleph _1$ countable Borel sets. The first case will be ruled out by
a theorem which says that Mathias forcing does not add a perfect set
of old reals. In the second case we will remember that really
$V=V[G_\a]$ for some $\a\in C$, and by the definition of $C$ we will
obtain a contradiction.

In order to prove that ${\cal Y}$ is uncountable, by fusion we will
build a perfect tree of $( N[G_0],P_\d /\name{G_0})$-generic conditions which
all decide $\D \cap N$ in different ways. This is much harder than it might
seem at first glance. The crucial lemma will be Lemma 4.7 below.

\sm

{\bf Definition 4.4}  (1) For $u\in [\delta ]^{<\o}$ and $p\in P_\delta$, let
$E(p,u)= \{ a\in ([\o ]^\o )^V : \exists q\leq _u p (q\Vdash
_{P_\d }\, a\in \D )\} $.

(2) Suppose $\bar x =\l \name{ x_\a} :\a \in u\r $ is such that every
$\name{ x_\a} $ is a $P_\a $-name for a finite subset of $\o $ with elements
larger than the members of the first coordinate of $p(\a )$. Then by
$p\cup \bar x$ we denote the condition $\bar p\in P_\d $ with $\bar
p(\a )=p(\a )$ for $\a \not\in u$, and first coordinate of $\bar p(\a ) =$
first coordinate of $p(\a )$, and second coordinate of $\bar p (\a ) =
($second coordinate of $p(\a ))\cup \name{ x_\a} $, for $\a \in u$.
Moreover, by $\bar x \cup p$ we denote the condition $\bar q \in P_\d
$ with $\bar q(\a )= p(\a )$ for $\a \not\in u$, first coordinate of
$\bar q(\a )=$ (first coordinate of $p(\a )$) $\cup \dot x_\a $ and
second coordinate of $\bar q(\a )=$ (second coordinate of $p(\a )$)
$\setminus (\max (\name{x_\a} ) +1)$ for $\a \in u$.

\sm

{\bf Lemma 4.5} {\it The ordering $\leq _u$ has the pure decision
property, that is, for $\tau $ a $P_\d $-name for a member of $\{ 0,1\} $
and $p\in P_\d $ there exists $q\leq _u p$ such that $q$ decides $\tau $.}

\sm

{\it Proof:} We prove it by induction on $\max (u)$. Let $\a _0=\max
(u)$ and $u_0=u\setminus \{ \a _0\} $. We may regard $\t $ as a $P_{\a
_0}$-name for a $P_{\d }/\Name{G_{\a _0}}$-name. Firstly, if $\a
_0=0$, then by the pure decision property of Mathias forcing (proved
in [B, 9.3.])  there exists $q(0)\in Q$, $q(0)\leq _{\{ 0\} }p(0)$,
deciding the disjunction ``$\exists q_1\in P_{\d }/\name{G_0}(q_1\leq
p\upharpoonright [1,\d ) \wedge q_1\forces _{1\d }\, \t =0)\vee
\exists q_1\in P_{\d }/\name{G_0} (q_1\leq p\upharpoonright [1,\d )
\wedge q_1 \forces _{1\d }\, \t =1 )$''. By the maximum principle of
forcing we
may find $q_1$ such that $q(0)\extend q_1 \leq _{\{0\}} p$ and
$q(0)\extend q_1 $ decides $ \t $. 

For the inductive step, as in the case $\a _0 =0$ we know that for
some $q_1 \in P_{\d }/\Name{G_{\a _0}}$, $q_1\leq _{\{ \a _0\}
}p\upharpoonright [\a _0, \d )$, $p\upharpoonright \a _0 \forces
_{P_{\a _0}}\, ``q_1
\hbox{ decides } \t $''; moreover, by induction hypothesis there exists
$q_0\leq _{u_0} p$, $q_0\in P_{\a _0 }$, which decides whether for
such $q_1$, $q_1\forces \, \t =0 $ or $q_1 \forces \t =1 $. Then
$q_0\extend q_1$ is as desired. \hfill $\qed $

\sm

{\bf Lemma 4.6.} {\it Let} $p\in P_\d $, $u\in [\dom (p)]^{<\o }$,
$n\in \o $ {\it and $\bar x=\l \name{ x_\a} :\a \in u\r $ such that
$\name{ x_\a} $ is a $P_\a $-name for the first $n$ members of the
infinite part of $p(\a )$. Suppose also that for no $q\leq p$,
$E(q,u)$ is a filter.

Then for $i\in \{ 0,1 \} $ there exist $q_i\leq _u p$ and disjoint
$a_i \in \inf $ such that} $ q_i \cup \bar x  \Vdash ``a_i \in  \D
\hbox{''}.$ 

\sm

{\it Proof:}  First note that if $q\leq  p$, for every $k\in \o $ we may
find a disjoint sequence $\l a_i:i<k\r $ of members of $\inf $ and $\l q_i:
i<k\r $ such that $q_i \leq _u q$ and $q_i\Vdash ``a_i \in \D \hbox{''}$.
In fact, since $E(q,u)$ is not a filter there exist $a_0', a_1'\in E(q,u)$
such that $a_0'\cap a_1'\not\in E(q,u)$. Let $q_i'\leq _u q$ force ``$a_i'
\in  \D $''. By the pure decision property of $\leq _u$, as proved in
Lemma 4.5., there exists 
$q_0\leq _u q_0'$ deciding whether $a_0:=a_0'\setminus a_1'$ or $a_0'\cap
a_1'$ belongs to $\D$. But then clearly $q_0\Vdash ``a_0 \in  \D
\hbox{''}$. Hence we may let $q_1=q_1'$, $a_1=a_1'$. 
Now proceeding by induction we easily construct $\l a_i:i<k\r $ and
$\l q_i:i<k\r $ as desired.

For $\a \in u$ let $\l \Name{ y^i_\a} :i<2^n\r $ be an enumeration (of names)
of all the subsets of (the denotation) of $\name{ x_\a} $, and let 
$\l \bar y_i:i<n^*\r $ enumerate all $\bar y_\sigma =\l \Name{ y ^{\sigma (\a
)}_\a} :\a \in u\r $, where $\sigma \in {^u(2^n)}$. Now using the observation
above we easily construct $q_\tau $ and $a_\tau \in \inf $, for every
$\tau \in {^{\leq n^*}(n^*+1)}$, such that the following requirements hold:

\item{(1)}{$q_\emptyset =p$, $a_\emptyset =\o ,$}

\item{(2)}{$\l a_{\tau \hat {\; } \l i\r } : i< n^*+1\r $ is a partition of 
$\o $,}

\item{(3)}{$\tau \sub \sigma \Rightarrow q_\tau \geq _u q_\sigma $,}

\item{(4)}{$|\tau |>0\Rightarrow \bar y_{|\tau |-1} \cup
q_\tau \Vdash
``a_\tau \in  \D \hbox{''}$.}

Now choose $q_0\leq _up$ such that for every $i<n^*$ and $\tau \in {^{<n^*}
(n^*+1)}$, $\bar y_i \cup q_0 $ decides for which $j$, $a_{\tau \hat {\; }
\l j\r }$ belongs to $\D$. For this we use again the pure decision property of
$\leq _u$. Then clearly we may find $\tau _1\in {^{n^* }(n^* +1)}$ such that,
letting $a_1:= \bigcup \{ A_{\tau _1 |j }: 1\leq j\leq n^*\} $, $a_0:=\o 
\setminus a_1 $ and $q_1 :=q_{\tau _1}$, the conclusion of the Lemma holds.
\hfill $\square $

\sm

The following lemma shows that the assumption of Lemma 4.6 holds.
As always, we implicitly regard $P_\d /\name{G_0}$ as a countable
support iteration of Mathias forcing.

\sm

{\bf Lemma 4.7.}  {\it In $V[G_0]$, for no $q\in
P_\d /\name{G_0}$ with $q\leq p^*\restrict [1,\d )$,  and for no} 
$u\in [\dom (q)]^{<\o }$ is it true that $E(q,u)$ {\it is
a filter.}

\sm

{\it Proof:}  Suppose by way of contradiction that for some $q\leq p^*\restrict
[1,\d )$
and $u\in [\dom (q)]^\o $, $E(q,u)$ is a filter. By the pure decision property
of $\leq _u$, then $E(q,u)$ is an ultrafilter. By the transitivity of
the ordering $\leq _u$ we have that for every $q'\leq _uq$,
$E(q',u)\sub E(q,u)$ and hence $E(q',u)$ is a filter. 
By the pure decision property again, we obtain $E(q',u)= E(q,u)$.
This fact will be used several times in the sequel.

In $V$ let $\E, \name{ q}$ be $Q_0$-names for $E(q,u), q$. 
Without loss of generality 
we may assume that the above properties of $E(q,u),q $ are forced by
$p^*(0)$ to hold for
$\E,\name{ q}$. Moreover we may certainly assume $\E, \name{ q}
\in N$.

Let $G_0=G'_0\ast G''_0$ be the decomposition of $G_0$ according to
the decomposition of Mathias forcing
$Q_0=Q_0'\ast \name{Q_0''}$. Let $p^*(0)=( u^{p^*},a^{p^*} ) $. 
In $V[G'_0]$ we can define:

$$D _1= \{ a\in \inf : \exists a' \in G'_0 ( u^{p^*} , 
a') \Vdash ``a\in \E\hbox{''}\} .$$

By hypothesis and as $Q(G'_0)$ has the pure decision property (see [JSh]), 
we conclude
that $D _1$ is an ultrafilter. Working in $V[G_0']$, 
we distinguish two cases according to whether
$G'_0$ is a projection of $D _1$ or not. In both cases we derive a 
contradiction:

\sm

{\it Case 1:}  $G'_0 \leq _{RK}D _1$. 

\sm

Let $f\in \fun $ witness this. As $Q_0'$ is $\sigma $-closed and hence
does not add new reals, $f\in V$. As $N':= N[G'_0]\prec ( H(\chi
)^{V[G'_0]}, \in ) $ (see [Sh{\it b}, 2.11., p.88]) and $D _1 \in
N'$, we may assume $f\in N'$, and hence $f\in N$ by properness. As
$G'_0 \cap N$ is countable, there exists $a\in G'_0$ such that
$G'_0\cap N=\{ b\in N:a\sub ^*b\}$. 

We work in $V[G_0']$. By Case 1 there exists $b\in D _1$ such
that $f[b]\sub a$. Let $x\in Q(G_0')$ with $u^{p^*}$ as its first
coordinate be such that 
$$x\forces _{ Q(G_0')}\, b\in \E.\eqno(1)$$

\no Note that $x$ is trivially $(N[G_0],Q(G_0'))$-generic, since
$Q(G_0')$ is {\it ccc}. By Lemma 3.10 there exists a $Q(G_0')$-name
$\name{ q_1}$ for a $(N[\name{ G_0}],P_\delta /\name{ G_0})$-generic condition,
such that $p^*(0)$ $\forces _{Q(G_0')}$ $\name{ q_1}\leq_u\name{ q}$. By the
remark at the beginning of this proof, we have

$$p^*(0)\forces _{Q(G_0')}\, E(\name{q},u)=E(\name{ q_1},u).\eqno(2)$$

We conclude that $x\ast \dot q_1$ is a $(N[G_0'],Q(G_0')\ast (P_\delta
/\name{ G_0}))$-generic condition below $p^*$. By (1) and (2), there is a
$Q(G_0')$-name $\name{ q_2}$ such that $p^*(0)\forces\, \name{ q_2}\leq
_u\name{ q_1}$ and $x\ast \name{ q_2}\forces \, b\in \D$. Clearly, 
$x\ast\name{ q_2}$ is $(N[G_0'],Q(G_0')\ast (P_\delta
/\name{ G_0}))$-generic. Let $G_1$ be
$Q(G_0')\ast (P_\delta /\name{ G_0})$-generic over
$V[G_0']$ such that $x\ast \name{ q_2}\in G_1$.
We conclude that $b\in \D [G'_0\ast G_1]$. 

Clearly we have that $f_\ast (\D [G'_0 \ast G_1])\ne G'_0 $, since
otherwise $\D [G'_0\ast G_1]$ could be computed from $f$ and $G'_0$ in
$V[G'_0]$. For this we use that $\D [G'_0 \ast G_1] $ is Ramsey. Hence this
inequality holds in $N'[G_1]$. Therefore there exists $a_1\in N\cap
G_0'$ such that $f^{-1}[a_1]\not\in\D [G'_0\ast G_1]$. Let
$b_1=b\setminus f^{-1}[a_1]$. So $b_1\in \D [G'_0\ast G_1]$. We obtain
that $f[b_1]\cap a_1=\emptyset$, $f[b_1]\sub a$, and $a\sub ^*a_1$.
Hence $f[b_1]$ is finite. But then $f[b_1]\not\in G'_0$, a
contradiction.

\sm

{\it Case 2:}  $G'_0\not\leq _{RK} D _1$.   

\sm

In $V$ let $\D _1$ be a $Q_0'$-name for $D _1$, and 
let $\name{ G_0'}$ be the canonical name
for the $Q_0'$-generic filter. Then by hypothesis there exists $t_0 
\in [a^{p^*}]^\o $ such that:

$$t_0 \Vdash _{Q_0'} ``{\name{ G_0'}} \not\leq _{RK} \D _1 \hbox{''}.$$ 

\no We may certainly assume $ \D _1, t_0 \in  N$. 

In $V$ let $g$ be $Q'$-generic over $N$ such that $t_0\in g$, where
$Q$ is Mathias forcing and $Q=Q'\ast \name{Q''}$ its canonical
decomposition. In $N[g]$ let $d= \D _1[g ].$ By elementarity we
conclude

$$N[g]\models d\hbox{ is an u.f.}, g \hbox{ a Ramsey u.f. and }
g\not\leq _{RK} d .\eqno(3)$$

In [GSh] it was shown that for any ultrafilter $D$ on $\o $ there
exists a proper forcing $Q_{D}$ such that whenever
$G$ is a Ramsey ultrafilter with $G\not\leq _{RK} D$, then after forcing
with $Q_D$, $G$ still generates an ultrafilter but $D$ does not.
Moreover $Q_D$ is $\fun $-bounding. Hence by Lemma 1.1, every such $G$ 
generates a Ramsey ultrafilter in every $Q_D$-generic extension. 

\sm

{\bf Definition 4.8.} Conditions in $Q_D$ are $f=\l h,E;
E_0,E_1,\dots \r $ where $h:\o \rightarrow \{ -1,1\} $, and the sets
$E,E_0,
E_1,\dots $ belong to the ideal dual to $D$ and partition $\o $. 

The ordering is defined as follows: $\l h,E;E_0,E_1,\dots \r \leq \l
h',E';E_0', E_1',\dots \r $ if and only if 

$E\supseteq E',$

$E,E_0,E_1,\dots $ is a coarser partition than $E',E_0',E_1',\dots ,$

$h\restrict E' =h'\restrict E',$ 

for all $i$: $h\restrict E_i' \in \{ h'\restrict E_i' ,-h'\restrict E_i' \} $.

\sm

A $Q_D$-generic filter $G$ determines a generic real $s= \bigcup \{
h_f: f\in G\} $.

By standard arguments one proves 
that whenever $s\in {^\o \{ -1,1\} }$ is $Q_{D}$-generic,
$f$ belongs to the generic filter which $s$ generates, and $s_f$ is defined by:

$$s_f(n)= \cases{s(n) &$n\in E^f $\cr
                -s(n)&$n\not\in E^f ,$ \cr}\eqno(2)$$

\noindent then $s_f$ is $Q_{D}$-generic as well and $f$ belongs to its
generic filter. Here $E^f$ is the second coordinate of $f$. Hence
especially $-s$, where $(-s)(n)=-s(n)$, is also
$Q_{D}$-generic.
 
\bigskip

In $N[g]$ we have the forcing $Q_{d}$.
In $V$, choose $s\in {^\o \{ -1,1\} }$ $Q _{d}$-generic 
over $N[g]$. By the properties of $Q_d$ and (3),
$g$ generates a Ramsey ultrafilter in $N[g][s]$.

Finally, in $V$ choose $t_1\sub t_0$ $Q(g)$-generic over
$N[g][s]$. Since every
infinite subset of $t_1$ is also $Q(g)$-generic 
and, as just noticed, $-s$ is also $Q _{d}$-generic,
without loss of generality we may assume that

$$V\models t_1 \Vdash _{Q'_0} ``s^{-1}(1)\in  \D _1 \hbox{''}.$$

\noindent Otherwise work with some $t_2\in [t_1]^\o $ and $-s$.
Hence, by the definition of $D _1$, and since $Q'_0$ does
not add reals, we may assume:

$$V\models ( u^{p^*}, t_1) \Vdash _{Q _0} ``s^{-1}(1)\in \E\hbox{''}.
\eqno(5)$$

\bigskip

{\bf Claim 1.}  {\it There exists a $Q $-name 
$\name{ q'} \in N[g,s, t_1] $ such that} 
$$N[g,s, t_1] \models ( u^{p^*}, t_1) \Vdash _{Q }
``\name{ q'} \in P_{\d }/\name{ G_0} \wedge \name{ q'} \leq _u \name{ q} \wedge
\name{ q'} 
\Vdash _{P_\d /\name{ G_0}} `\name{ r}
\alm s^{-1}(1)\hbox{'} \hbox{''}.\eqno(6)$$

\bigskip

{\it Proof:} Otherwise there exist $( u',t') \leq ( u^{p^*}, t_1)$ and
$\name{ q '}$ such that in $N[g,s,t_1]$, $\name{ q'}$ is a $Q$-name for a
condition in $P_\d /\name{ G_0}$, and

$$N[g,s, t_1]\models ( u',t') \Vdash _{Q } ``\name{ q'}\leq _u \name{ q} 
\wedge \name{ q'} \Vdash _{P_\d/\name{ G_0}} \hbox{`} \name{ r} \setminus s^{-1}(1)
 \hbox{ is infinite'} \hbox{''}.$$

\noindent Such $\name{ q'}$ exists by the existential 
completenes of forcing and the pure decision property of $\leq _u$.

By Lemma 3.10, in $V$ there exists $\bar q \in P_\d$ such that $\bar q
\leq _u ( u',t')\ast \name{ q'} $ and $\bar q$ is $( N[g,s,t_1], P_\d
)$-generic.  Since by the observation at the very beginning of the
present proof we know that

$$p^*(0) \Vdash _{Q } `` \E = E(\bar q\restrict [1,\d ),u) \hbox{''},$$

\noindent by (5) and the definition of $\E$, there exists $\bar {\bar q} \in
P_\d $ such that $\bar {\bar q}
\leq _u \bar q $ and $\bar {\bar q} \Vdash _{P_\d }``\name{ r} \alm
s^{-1}(1)$''. Now choose $G$ $P_\d $-generic over $V$ such that $\bar {\bar
 q}\in G$.
Then clearly $V[G]\models \name{ r}[G]\alm s^{-1} (1)$ 
and $N[g,s,t_1][G] \models |\name{ r} [G]
\setminus s^{-1}(1)| =\o $. But $\name{ r}[G]$ is the same 
real in
both models, a contradiction. \hfill$\square$

\sm

Let us abbreviate the formula ``$\dots $'' in (6) by $\phi (\name{ q'}, s)$.  

Since $t_1$ is $Q(g)$-generic and $g$ generates a Ramsey 
ultrafilter in $N[g][s]$, there exists $( u',t' ) \in Q(g)$ 
such
that $u'\sub t_1 \sub t'$ and

$$N[g,s]\models ( u', t') \Vdash _{Q(g)} ``( u^{p^*},\name{ t}) 
\Vdash _{Q }  \hbox{`} \phi (\name{ q'}, s) \hbox{' ''},\eqno(7)$$

\noindent where $\name{ t}$ is the canonical name for the generic real added by
$Q (g)$, and in the formula $\phi (s)$, $\name{ q'}$ is now a
$Q(g)$-name for the above $\name{ q'}$. 

Since $s$ is $Q _{d}$-generic over $N[g]$ and $( u',t') \in 
N[g]$, there exists $f
\in Q _{d}$ such that $f$ belongs to the 
$Q _{d}$-generic 
filter induced 
by $s$, and in $N[g]$ the following holds:

$$f\Vdash _{Q _{d}}``N[g][\name{ s}]\models (( u',t') \Vdash _{Q(g)}
`( u^p,\name{ t}) 
\Vdash _{Q  }\, \phi (\name{ q'},\name{ s}) \hbox{'})\hbox{''},$$

\noindent where $\name{ s}$ is the canonical $Q _{d}$-name for the 
$Q _{d}$-generic real and in $\phi (\name{ q'},\name{ s})$, $\name{
q'}$ denotes now a $Q_{d}\ast Q(g)$-name for the $\name{ q'}$ in (7).
By the definition of $Q_{d} $ we have $\o \setminus E^f \in d$.

\sm

{\bf Claim 2:}   $V\models t_1 \Vdash _{Q' } ``\o \setminus 
E^f \in  \D _1 \hbox{''}.$

\sm

{\it Proof:}  As $g$ is $Q'$-generic over $N$, $\o \setminus 
E^f
\in d =  \D _1 [g]$, and $Q' $ does not add reals, there exists 
$u\in g $ such that 

$$N\models  u\Vdash _{Q' } ``\o \setminus E^f \in  \D _1 \hbox{''}.$$

\no By elementarity we conclude that this is true in $V$. But clearly we have 
$t_1\alm u$. \hfill$\square$

\sm

Let $s_f$ be defined as in definition 4.8.  By the remarks after 4.8,
$s_f$ is $Q _{d}$-generic over $N[g]$, and clearly $f$ belongs to the
generic filter determined by $s_f$. Hence (5) holds if $s$ is replaced
by $s_f$. Clearly $N[g][s]=N[g][s_f]$, and hence $t_1$ is
$Q(g)$-generic over $N[g][s_f]$, and consequently
$N[g][s][t_1]=N[g][s_f][t_1] =:N^*$.

Let $G^*$ be $Q$-generic over $V$, containing a $( N^*, 
Q ) $-generic condition below $( u^{p^*}, t_1) $. Then by Claim 2, $\o 
\setminus E^f\in \E [G^*]$. But also $s^{-1}(1), s_f^{-1}(1)\in \E[
G_0^*]$. In fact, in $N^*[G^*]$ we have $q_1:=\name{ q'}[s][t_1][G^*]$
and $q_2: =\name{ q'}[s_f][t_1][G^*] $ with the property that $P_\d
/\name{ G_0}\models q_1, q_2 \leq _u q$, and $q_1 \Vdash _{P_\d
/\name{G_0} }\,
\name{ r} \alm s^{-1}(1)$, and $q_2 \Vdash _{P_\d /\name{G_0} } \, 
\name{ r} \alm
s_f^{-1}(1)$.  Otherwise, as in the proof of Claim 1, in $V[G^*]$ we
could find $(N[G^*],P_\d /\name{ G_0})$-generic conditions $\bar
q_1\leq _u q_1$ and $\bar q_2\leq _u q_2$ forcing the opposite. By
choosing filters which are $P_\d /\name{ G_0}$-generic over $V$ and
contain $\bar q_1,\bar q_2$ respectively, we obtain a contradiction.
Consequently, $s^{-1}(1)$, $s_f^{-1}(1) $, and $\o
\setminus E^f$ belong to $\E [G^*]$. But $s^{-1}(1), s_f^{-1}(1)$ 
are complementary on $\o
\setminus E^f$, and hence $\E [G^*]$ is not a filter, a
contradiction. \hfill $\square $

\sm

Using 4.6, 4.7 and [B, Lemma 7.3.], by standard arguments on
proper forcing we obtain the
following Corollary.

\sm

{\bf Corollary 4.9.}  {\it In $V[G_0]$, there exist $\l q_s:s\in 
{^{<\o } 2} \r , \; 
\l a_s :s\in {^{<\o }2}\r $ such that the following hold: }

\item{(1)}{if $s\sub t$, then $ a_s\supseteq a_t$ and 
$a_{s\hat {\; }\l 0\r }
\cap a_{s\hat {\; } \l 1\r } = \emptyset ,$}

\item{(2)}{{\it if $f\in {^\o 2}$, then $\l q_{f\restrict n}:n<\o \r $ is a
descending chain in $P_{\d }/\name{G_0}$ which has a lower bound 
$q_f$ such that:}} 

\itemitem{$\cdot $}{$q_f$ {\it is $(N[G_0] ,P_{\d }/\name{G_0} )$-generic,}}

\itemitem{$\cdot $}{$q_f \Vdash \forall n (a_{f\restrict n}\in  \D ) $,}

\itemitem{$\cdot $}{$q_f $ {\it decides $ \D \cap N$.}}

\sm

{\bf Corollary 4.10.}  {\it $\Y $ is uncountable.}

\sm

\relax From Lemma 4.3 and Corollary 4.10 we conclude that $\Y $ is an
uncountable $\Sigma ^1_2$ set in $V[G_0]$ which is a subset of $V$. By
well-known results from descriptive set theory, $\Y
$ is the union of $\o _1 $ Borel sets, say $\l B_\a :\a < \o _1 \r $,
and this decomposition is absolute for models computing $\o _1 $
correct (see [J, Theorem 95, p.520, its proof on p.526 using the
Shoenfield tree, and Lemma 40.8, p.525, where its absoluteness is
proved]). If one of the $B_\a $ is uncountable it contains a perfect subset
(see [J, Theorem 94, p.507]). This case will be ruled out by the next Lemma
4.11.

Otherwise, each $B_\a $ is countable. Now $\Y $ and hence $\l B_\a :\a
< \o _1 \r $ is coded by a real $x$. We may assume that $x$ also codes
$\l q_s :s\in {^{<\o }2} \r $ and $\l a_s : s\in {^{<\o }2} \r $ from
4.9. Now remember that $V$ here is really $V[G_\a ]$ where $\a \in
C$ ($C$ coming from 4.1), and hence $V[G_0] =V[G_{\a +1}]$. Clearly there
exists $\b < \a $ such that $x\in V[G_\b ,r_\a ]$. Then also $\l B_\a
:\a <\o _1 \r $, $\l q_s :s\in {^{<\o }2}\r , \l a_s :s\in {^{<\o
}2}\r \in V[G_\b ,r_\a ]$, and hence, as each $B_\a $ is countable, $\Y
\sub V[G_\b ,r_\a ]$. But from this we conclude $\P ^{V[G_\b ,r _\a ]}
=\P ^{V[G_{\a +1}]}$, as a new real in $V[G_{\a +1}]\setminus V[G_\b
,r_\a ] $ would give a new branch through $\l a_s :s\in {^{<\o }2} \r
$ and hence a new member in $\Y $. But $\a \in C$, and hence $(\ast )$
in 4.1 fails for it, a contradiction.

\sm

Therefore, 
in order to finish the proof of Proposition 2.3 it suffices to prove the
following Lemma:

\sm

{\bf Lemma 4.11.}  {\it Suppose $q\in Q$, where $Q$ is Mathias-forcing,
$\t $ is a $Q $-name, and

$$q \Vdash _Q ``\t \sub \tree \hbox{ is a perfect tree''}.$$

\noindent Then $q\Vdash _Q ``[\t ]\not\subseteq V\hbox{'' }$.}

\sm

{\it Proof:} By applying the pure decision property of $Q $ 
repeatedly, without loss
of generality we may assume that if $q=( s,a) $, then for every $t\in [a]^{
<\o }$ and $n\in \o $ there exists $m\in \o $ such that $( s\cup t,a
\setminus m) $ decides the value of $\t \cap {^{<n}2}$. Hence if we let

$$T_t=\{ \nu \in \tree : \exists n(( s\cup t, a\setminus n) \Vdash _Q
``\nu \in \t \hbox{''})\} ,$$

\noindent then $T_t$ is a tree with no finite branches.

We shall define a $Q $-name $\eta $ for a real in $[\tau ]\setminus
V$. To this end, for every
$t\in [b]^{<\o }$, we construct   
$b\in [a]^\o $, $\eta _t \in T_t$ and $n(t)\in \o $ such that the
following hold:

\smallskip

\item{(1)} $( s\cup t,b\setminus (\max (t)+1 )) \Vdash _Q ``\eta _t
\in \t \hbox{''}$;

\smallskip

\item{(2)}if $T_t\cap [\eta _t]$ has infinitely many branches, hence
by compactness a nonisolated one, and $x_t$ is the lexicographically least
such one, then 
for every $m\in b\setminus (\max (t)+1)$, $\eta _{t\cup \{ m\} }$ is
not an initial segment of $x_t$, but $\eta _{t\cup \{ m\} }\restrict n(t\cup \{
m\} )$ is; moreover, $\lim _{m\rightarrow \o }n(t\cup \{ m\} )=\infty $;

\smallskip

\item{(3)}if $T_t\cap [\eta _t]$ has finitely many branches, then
for every $m\in b\setminus (\max (t)+1)$,

\itemitem{$\cdot $}if $T_{t\cup \{ m\} }$ has a member extending $\eta _t$
which does not belong to $T_t$, then $\eta _{t\cup \{ m\} }$ is like that, 
say among the shortest the lexicographically least one;

\itemitem{$\cdot$}if $T_{t\cup \{ m\} }$ has no such member, then
$\eta_{t\cup \{ m\} }=\eta _t$.

\smallskip

The construction of $b$, $\l \eta _t:t\in [b]^{<\o }\r $ and $\l n(t):
t\in [b]^{<\o }\r $ is by fusion:
Suppose that an initial segment of $b$, say $t$, has been fixed and for some
$b'\in [a\setminus t]^\o $, for every $t'\in {\cal P}(t )$ and $m\in b'$,
$\eta _{t'}, n(t')$ and $\eta _{t'\cup \{ m\} }, n(t'\cup \{ m\} )$ 
have been defined such that (1),
(2), (3) hold for $\eta _{t'}, \eta _{t'\cup \{ m\} },
n(t'), n(t'\cup \{ m \} )$ and $b'$. 
Now the least element of $b'$, say
$k$, is put into $b$. Then successively for each $t'\in {\cal P}(t )$,
first count how many branches $T_{t'\cup \{
k\} }\cap [\eta _{t' \cup \{ k\} }]$ has, and then accordingly define
$\eta _{t' \cup \{ k\} \cup \{ m\}}$ and maybe $n(t'\cup \{ k\} \cup
\{ m\} )$ (if
we are in case (2)) for $m\in b'$, all the time shrinking
$b'$ to make sure that in the end, for some $b''\in [b']^\o $, for
every $t'\in
{\cal P}(t\cup \{ k\} )$ (1), (2) and (3) hold for $\eta _{t'}$ and
$b''$. The construction is totally straightforward, so we leave
the rest to the reader.

We define a $Q$-name as follows:

$$\eta = \bigcup \{ \eta _t:t\in [b]^{<\o } \wedge 
(( s\cup t , b\setminus (\max (t)+1) ) \in \name{ G} )\} . $$

Here $\name{ G}$ is the canonical name for the $Q $-generic filter. By 
construction
we conclude:

$$( s, b) \Vdash _Q ``\eta \in [\t ]\cup \tau \hbox{''}.$$

Suppose now that some 
$( s\cup t, b^*) \leq ( s,b) $ forces that $\eta $ belongs to 
$V$, so, without loss of generality, there exists $\eta ^* \in V$ such that 

$$( s\cup t,b^*) \Vdash _Q `` \eta =\eta ^* \hbox{''}.$$

\noindent \relax From this we will derive a contradiction. Then the Lemma
will be proved. Clearly we have $\eta ^* \in \fun \cup \fin $. We 
distinguish the following cases:

\smallskip

{\it Case 1:}  $T_t\cup [\eta _t]$ has infinitely many branches.

{\it Subcase 1a:}  $\eta ^* =x_t$. By construction, if $m\in b^*$ then
$( s\cup t\cup \{ m\} , b^*\setminus (m+1))$ $ \Vdash _Q$     
$``\eta _{t\cup \{ m\} }
\sub  \eta \hbox{''}$ and $\eta _{t\cup \{ m\} } \not\subseteq x_t$,
a contradiction.

{\it Subcase 1b:}  $\eta ^*\restrict n \ne x_t\restrict 
n $ for some $n$. If $m\in b^*$
with $n(t\cup \{ m\} )\geq
n$, then by construction $( s\cup t\cup \{ m\} , b^* \setminus
(m+1)) \Vdash _Q \,
``\eta \restrict n(t\cup \{ m\} ) =x_t\restrict 
n(t\cup \{ m\} )\hbox{''}$, a contradiction.

\smallskip

{\it Case 2:}  $T_t\cap [\eta _t]$ has only finitely many branches. 

{\it Subcase 2a:}  $\eta ^* \in [T_t]\cup T_t$. Since $\t $ is forced to be a
perfect tree there exists $u\in [b^*]^{<\o }$ such that $T_{t\cup u}$ has
a member above $\eta _t$ which is not in $T_t$. But then by construction
$( s\cup t\cup u, b^*\setminus (\max (u)+1)) \Vdash _Q \,
\hbox{``}  \eta 
\not\in [T_t]\cup T_t \hbox{''}$, a contradiction.

{\it Subcase 2b:} $\eta ^*\restrict n \not\in T_t$ for some $n$. By
construction of $T_t$, there exists $m$ such that $( s\cup t,
b^*\setminus m )
\Vdash _Q \,
``\t \cap {^{\leq n}2} = T_t \cap {^{\leq n}2}$''. But
$( s\cup t, b^*\setminus m ) \Vdash _Q \, 
\hbox{``} \eta \restrict n\in \t \hbox{''}$,
a contradiction. 

\hfill $\square \square $

\bigskip

{\hgross 5. Proof of Proposition 2.4}

\bigskip

The proof will use several ideas from the proof of Proposition 2.3.
Suppose that Proposition 2.4 is false, that is, 
there exist $Q$-names
$\D $ and $ \name{ r}$, and $p\in
Q $ such that $p$ forces that $\name{ r}$ induces a Ramsey ultrafilter 
$\D $ on $(\inf )^V$ which is not RK-equivalent to $\name{ G'}$ by any
$f\in \fun\cap V$.

First note that a $\s $-centered forcing $P$ does not add
such $ \D $. In fact, 
since $V\models CH$, 
such $ \D $ is forced to be generated by a $\alm
$-descending chain $\l \name{a_\a} :\a < \o _1\r $ of members of $(\inf
)^V$. For every $\a <\o _1$, choose $p_\a \in P$ and $a_\a \in (\inf
)^V$ such that $p_\a \Vdash _P \; \name{ a_\a} =a_\a $. Since $P$ is $\s
$-centered, there exists $X\in [\o _1]^{\o _1}$ such that $p_\a ,p_\beta $
are compatible whenever $\a , \beta \in X$. By the {\it ccc} of $P$, there
exists a $P$-generic filter $G$ which contains $p_\a $ for uncountably
many $\a \in X$. Then clearly $ \D [G]\in V$, as $\D [G]$
is generated by $\l a_\a :\a \in X\r $. The argument shows that no
condition in $P$ forces that $\D $ does not belong to $V$. 

Since $Q(\name{ G'})$ is forced to be $\s$-centered, by what we just
proved we may assume that $\D$ is a $Q'$-name. 
As usual, we write $p= ( u^p, a^p ) $. For $t\in Q'$ we define
$$D _t = \{ a\in (\inf )^V : t\Vdash _{Q'} ``a\in {\D} 
\hbox{''} \}. $$

The following claim follows immediately from the definitions:  
 
\sm

{\bf Claim 1.} {\it For all $t\in Q' $ with $t\leq a^p$, we have that
$a\in D _t$ if and only if $( u^p , t )$ $\Vdash_{Q }$ $``\name{ r}
\alm a \hbox{''} .$}

\sm

{\bf Claim 2.}  {\it Suppose that $( N, \in ) $ is a countable model of 
$ZF^-$ such that $\name{ r}, p\in N $, and $\name{ r}$ is hereditarily
countable in $N$. Then for every $a \in 
\inf \cap N$ and $t\in Q'\cap N $ with $t\leq  a^p$, it is true that
$( u^p , t) \Vdash _{Q }\, ``\name{ r} \alm a \hbox{''}$ implies that
$N\models ( u^p , t ) \Vdash _{Q }\, ``\name{ r} \alm a\hbox{''}.$}

\sm

{\it Proof of Claim 2:}  Otherwise there exists $q\in N\cap Q$ such
that $q\leq ( u^p , t )$ and
$N\models \, q\Vdash _{Q  }\, ``\name{ r}\cap (\o \setminus a) \hbox{ is
infinite}\hbox{''})$.  

By Lemma 1.2, there exists $q'\in Q$ such that $q' \leq q$ and $q'$ is
$(N, Q )$-generic. Let $G$ be $Q$-generic over $V$, containing $q'$.
Then by assumption $\name{ r}[G] \alm a$. On the other hand, $N[G]
\models |\name{ r}[G ]
\cap (\o \setminus a)|=\o $. As $\name{ r}[G]$ is
the same real in $V[G]$ and $N[G]$ we have a contradiction.\hfill $\square $

\def \d {\delta }
\sm

By assumption, and since $Q'$ does not add reals, we conclude:
$$a^p \Vdash _{Q'} \,``\D\hbox{ and }\name{ G'} \hbox{ are Ramsey
ultrafilters which are not }RK\hbox{-equivalent.''}$$ 

Choose a countable elementary substructure $( N,\in )\prec ( H(\chi )
, \in ) $ where $\chi $ is a large enough regular cardinal, such that
$\D , \name{ r} , p \in N$. 

In $V$, let $g$ be $Q'$-generic over $N$ such that $a^p \in g$. In
$N[g]$, let $d= \D [g].$ By elementarity we conclude 

$$N[g]\models ``g \hbox{ and }d\hbox{ are Ramsey
ultrafilters which are not }RK\hbox{-equivalent.''}\eqno(1)$$ 

In $V$, choose $s\in {^\o \{ -1 ,1\} }$ $Q_{d}$-generic over $N[g]$,
where $Q_{d}$ is the forcing from 4.8, defined in $N[g]$ from the
ultrafilter $d$. From (1), Lemma 1.1, and [GSh] we conclude that $g$
generates a Ramsey ultrafilter in $N[g][s]$.

Finally, in $V$ choose $t_1\leq a^p $ $Q(g)$-generic over
$N[g][s]$. Since every infinite subset of $t_1$ is also $Q(g)$-generic
and $-s $ is also $Q _{d}$-generic,
without loss of generality we may assume that $s^{-1}(1)\in D _{t_1}.$

By Claims 1 and 2 we conclude:
$$N[g][s][t_1] \models ( u^p ,t_1) \Vdash _{Q }\, ``\name{ r}{\alm
s^{-1}(1)}\hbox{''}.\eqno(2) $$

Since $g$ generates a Ramsey 
ultrafilter in $N[g][s]$, by the remark preceding Lemma 1.2 we
conclude that $t_1$ is $Q(g)$-generic over $N[g][s]$. Since
$Q(g)^{N[g]}$ is dense in $Q(g)^{N[g][s]}$, there exists $( u',t' )
\in Q(g)^{N[g]} $ 
such that $u'\sub t_1 \sub t'$ and
$$N[g,s]\models ( u', t') \Vdash _{Q(g)}\, ``( u^p ,\name{ t}) \Vdash
_Q\, `\name{ r} \sub ^* {s^{-1}(1)}\hbox{' ''}.\eqno(3) $$

\noindent Here $\name{ t}$ is the canonical name for the generic real added by
$Q(g)$.

Since $s$ is $Q _{d}$-generic over $N[g]$ and all the parameters
in the formula $`` \dots $'' of (3) belong to $N[g]$, there exists $f
\in Q _{d}$ such that $f$ belongs to the $Q _{d}$-generic 
filter induced 
by $s$, and in $N[g]$ the following holds:
$$f\Vdash _{Q _{d}}\, ``N[g][\name{ s}]\models [( u',t') \Vdash _{Q(g)}\,
`( u^p ,\name{ t}) \Vdash _{Q }\, `\name{ r} \sub ^* \dot s^{-1}(1)\hbox{''}]
\hbox{''}.$$

\noindent Here $\name{ s}$ is the canonical $Q _{d}$-name for the 
$Q _{d}$-generic real. 
By definition of $Q_d$, $\o \setminus E^f \in d$.

\sm

{\bf Claim 3:}  {\it $V\models \o \setminus E^f \in D _{t_1}$.}

\sm

{\it Proof of Claim 3:}  As $g$ is $Q'$-generic over $N$, $\o
\setminus E^f 
\in d = \D _1 [g]$, and $Q'$ does not add reals, there exists 
$w\in g $ such that 
$$N\models  w\Vdash _{Q'
 }\,  ``\o \setminus E^f \in \D \hbox{''}.$$

By elementarity we conclude that this is true in $V$, so by definition of
$D _w$, $\o \setminus E^f \in D _w$. Clearly we have $t_1\leq w$, so
$\o \setminus E^f \in D _{t_1}$.\hfill $\square $

\sm

Let $s_f$ be defined as in the remark after 4.8.
Then $s_f$ is also $Q _{d}$-generic over $N[g]$, and clearly 
$f$ belongs 
to the generic filter determined by $s_f$. Hence (3) holds if $s$ is replaced
by $s_f$. 

Clearly $N[g][s]=N[g][s_f]$, and hence $t_1$ is $Q(g)$-generic over
$N[g][s_f]$, and consequently $N[g][s][t_1]=N[g][s_f][t_1]$.

\relax From (3) we conclude:
$$N[g][s_f][t_1] \models ( u^p ,t_1) \Vdash _{Q }\, ``\name{ r}\alm
{s_f^{-1}(1)} \hbox{''}. \eqno(4)$$

\relax From Claim 3 together with Claims 1 and 2 we conclude:
$$N[g][s_f][t_1] \models ( u^p ,t_1) \Vdash _{Q }\, ``\name{ r}\alm \o
{\setminus E^f} \hbox{''}\, .\eqno(7) $$
              
Since $s^{-1}(1), s_f^{-1}(1)$ are complementary on $\o \setminus
E^f $,
(2), (4) and (5) imply that $\name{ r}$ is forced to be finite, a
contradiction. \hfill $\square \square \square $

\vfill \eject
\centerline{{\hgross References}}

\item{[B]}{{\capit J.E.Baumgartner}, {\it Iterated forcing,} in: Surveys in set
theory, A.R.D. Mathias (ed.), London Math. Soc. Lect. Notes Ser. 8,
Cambridge Univ. Press, Cambridge (1983), 1--59}

\item{[BPS]}{{\capit B. Balcar, J. Pelant, P.Simon,} {\it The space of
ultrafilters on $N$ covered by nowhere dense sets,} Fund. Math. 110
(1980), 11--24}




\item{[GSh]}{{\capit M. Goldstern, S. Shelah,} {\it Ramsey ultrafilters 
and the
reaping number -- Con$({\frak r}<{\frak u})$}, Ann. Pure Appl. Logic 49
(1990), 121--142.}

\item{[J]}{{\capit T. Jech,} {\it Set theory,} Academic Press, New York (1978).} 


\item{[JSh]}{{\capit H. Judah, S. Shelah,} {\it $\Delta ^1_2$-sets  
of reals,} Ann.
Pure Appl. Logic 42 (1989), 207--233.}

\item{[JSp]}{{\capit H. Judah, O. Spinas,} {\it Large cardinals and 
projective sets,} submitted.}

\item{[M]}{{\capit A.R.D. Mathias,} {\it Happy families,} Ann. Math. Logic 12
(1977), 59--111.} 

\item{[R]}{{\capit M.E. Rudin,} {\it Partial orders on the types of $\beta N$,}
Trans. AMS 155 (1971), 353--362.}
 
\item{[Sh{\it b}]}{{\capit S. Shelah,} {\it Proper forcing,} Lecture Notes in Math., vol.
942, Springer (1982).}

\item{[ShSp]}{{\capit S. Shelah, O. Spinas,} {\it The distributivity numbers of
finite products of ${\cal P}(\o )$/fin,} in preparation}

\bye